\documentclass[11pt,twoside]{article}

\usepackage{amsmath}
\usepackage{amssymb}
\usepackage{amsthm}
\usepackage[dvips]{graphicx}

\newtheorem{thm}{Theorem}[section]

\newtheorem{prop}[thm]{Proposition}
\newtheorem{cor}[thm]{Corollary}
\newtheorem{rem}[thm]{Remark}
\newtheorem{ex}[thm]{Example}

\title{\textbf{\Large Some examples of global Poisson structures on $S^4$}}
\author{\textsc{Takayuki MORIYAMA
and
Takashi NITTA
}}
\date{}
\begin{document}
\maketitle
\thispagestyle{empty}

\footnote{\textit{2010 Mathematics Subject Classification.} Primary 53D17; Secondary 37F35.}
\footnote{\textit{Key Words and Phrases.} 
Poisson structure, twistor method, holomorphic foliation. }

\begin{quote}
\textbf{\fontsize{10pt}{15pt}\selectfont Abstract.}
\textrm{\fontsize{10pt}{15pt}\selectfont 
A Poisson structure is represented by a bivector whose Schouten bracket vanishes. 
We study a global Poisson structure on $S^4$ associated with a holomorphic Poisson structure on $\mathbb{CP}^3$. 
The space of the Poisson structures on $S^4$ is a real algebraic variety in the space of holomorphic Poisson structures on $\mathbb{CP}^3$. 
We generalize it to $\mathbb{HP}^n$ by using the twistor method. 
Furthermore, we provide examples of Poisson structures on $S^4$ 
associated with codimension one holomorphic foliations of degree 2 on $\mathbb{CP}^3$. 
}
\end{quote}

\section{Introduction}
Classical mechanics is described as symplectic geometry in Hamiltonian formalism. 
%Symplectic structures are important from the view point of the quantization. 
The classical mechanics is also written by Poisson algebra, and it is formulated to Poisson geometry. 
Poisson geometry is a natural generalization of symplectic geometry. 
It has been developed to many branches in mathematical physics, 
for example, dynamical system, quantization, non-commutative geometry, reduction theory and global analysis (cf. \cite{AG}, \cite{V}, \cite{W}). 
Especially, in the global analysis it is an important problem how many Poisson structures exist on a manifold. 
For example, there is no symplectic structure on $S^4$ but Poisson structures exist. 
Our purposes of this paper are to investigate how many Poisson structures exist on $S^4$ and to construct explicitly many examples of Poisson structures on $S^4$. 
%One question is how many Poisson structures exist on $S^4$. 

A global antisymmetric 2-vector is called a bivector. 
The Schouten bracket of two bivectors is defined by the convolution of the covariant derivative of the one and the other, further antisymmetrization of it. 
The bivector is called a Poisson structure if the Schouten bracket of the bivector with itself vanishes. 
Since the Schouten bracket is valued in the space of antisymmetric 3-vectors, any bivector is always a Poisson structure on 2-dimensional manifolds. 
In general, the equation for vanishing the Schouten bracket is written as a system of partial differential equations. 
It is very complicated to solve the equation even on 3-dimensional manifolds. 
Moreover, it is much harder to solve the equation globally on the manifold. 
Hence, on 3-dimensional manifolds, the existence of the Poisson structure is non-trivial and the explicit construction is also difficult.  

In the holomorphic category, holomorphic bivectors and holomorphic Schouten brackets are considered (c.f. \cite{H1}). 
When the Schouten bracket vanishes, the holomorphic bivector is called a holomorphic Poisson structure. 
On complex surfaces, any holomorphic bivector is always a holomorphic Poisson structure 
since the Schouten bracket is also a holomorphic antisymmetric 3-vector. 
Hitchin and Goto had deep results for deformations of the holomorphic Poisson structures \cite{G}~\cite{H2}. 
On the 3-dimensional complex projective space $\mathbb{CP}^3$, a holomorphic Poisson structure corresponds to a codimension one holomorphic foliation of degree 2. 
Cerveau and Lins Neto proved that the space of holomorphic foliations of degree 2 has six irreducible components~\cite{CLN}. 
By their classification, a holomorphic Poisson structure on $\mathbb{CP}^3$ is completely classified (c.f. \cite{LPT}). 

In this paper, we consider global Poisson structures on the 4-dimensional sphere $S^4$ by using the twistor method \cite{GN}~\cite{NT}. 
The projective space $\mathbb{CP}^3$ is the Hopf fibration over $S^4$. 
We study Poisson structures on $S^4$ associated with holomorphic Poisson structures on $\mathbb{CP}^3$. 
We provide a real structure on the space of holomorphic Poisson structures on $\mathbb{CP}^3$. 
The space of the Poisson structures on $S^4$ has the $\mathbb{R}^*$-action induced 
by the $\mathbb{C}^*$-action on the space of holomorphic Poisson structures on $\mathbb{CP}^3$. 
Then we show that 
\begin{thm}\label{s1t1}
The space of Poisson structures on $S^4$ induced by holomorphic Poisson structures on $\mathbb{CP}^{3}$ is 
the real form of the space of holomorphic Poisson structures on $\mathbb{CP}^3$ and 
it is a real algebraic variety. 
Moreover, the quotient of the space by the $\mathbb{R}^*$-action is a real projective algebraic variety. 
\end{thm}

The sphere $S^4$ is regarded as the $1$-dimensional quaternionic projective space $\mathbb{HP}^1$. 
We extend Theorem~\ref{s1t1} to $\mathbb{HP}^n$ by the twistor method. 
The space $\mathbb{HP}^n$ does not admit any symplectic structure but Poisson structures. 
The space of holomorphic Poisson structures on $\mathbb{CP}^{2n+1}$ is realized as an algebraic variety. 
We obtain 
\begin{thm}\label{s1t2}
The space of Poisson structures on $\mathbb{HP}^{n}$ induced by holomorphic Poisson structures on $\mathbb{CP}^{2n+1}$ is 
the real form of the space of holomorphic Poisson structures on $\mathbb{CP}^{2n+1}$ and 
it is a real algebraic variety. Moreover, the quotient of the space by the $\mathbb{R}^*$-action 
is a real projective algebraic variety. 
\end{thm}

Especially, on $\mathbb{CP}^3$, a holomorphic Poisson structure induces a codimension one holomorphic foliation of degree 2. 
However the inverse transformation has much ambiguity. 
Pym provided an explicit correspondence of a $1$-form $\alpha$ induced a holomorphic foliation of degree 2 
to a unimodular quadratic Poisson bracket on $\mathbb{C}^4$~\cite{P}. 
The Poisson bracket induces the Poisson structure $w$ on $\mathbb{CP}^3$ given by 
\begin{equation*}
w=\sum_{i,j=0}^{3}\frac{dz_i \wedge dz_j \wedge d\alpha}{dz_0\wedge dz_1\wedge dz_2\wedge dz_3}\frac{\partial }{\partial z_i}\wedge \frac{\partial }{\partial z_j}\label{s1e1}
\end{equation*}
on $\mathbb{C}^4\backslash \{0\}$. 
We define a real structure on the space of $1$-forms induced holomorphic foliation of degree $2$ on $\mathbb{CP}^3$ and 
prove that the above correspondence preserves the real structures : 
\begin{thm}\label{s1t3}
The Poisson structure $w$ is real if the form $\alpha$ is real. 
\end{thm}
By the transformation, a ``real" holomorphic foliation on $\mathbb{CP}^3$ induces a Poisson structure on $S^4$. 
We obtain examples of such Poisson structures associated with Cerveau and Lins Neto's classification of holomorphic foliations on $\mathbb{CP}^3$ (see Section~\ref{s6.3}). 

This paper is organized as follows. 
In Section 3, we study holomorphic Poisson structures on $\mathbb{CP}^3$ and see that the space of such structures is an algebraic variety. 
In Section 4, we consider Poisson structures on $S^4$ associated with holomorphic Poisson structures on $\mathbb{CP}^3$ 
and show Theorem~\ref{s1t1} (Theorem~\ref{s4.5t1} and Corollary~\ref{s4.5c1}). 
In Section 5, we extend the argument to the higher dimensional case as Theorem~\ref{s1t2} (see Proposition~\ref{s5.1p1} and Theorem~\ref{s5.2t1}). 
In the last section, we show Theorem~\ref{s1t3} (Theorem~\ref{s6t1}) and provide examples of Poisson structures on $S^4$ related to holomorphic foliations on $\mathbb{CP}^3$.

\section{Preliminaries}
In this section, we explain the definition of Poisson structure. 
For much of this material, we refer to~\cite{H2} and \cite{V}. 
Let $M$ be a differentiable manifold of dimension $n$. 
A global section of $\wedge^q TM$ is called a {\it q-vector} on $M$. 
We write a 2-vector $u$ on $M$ by 
$u=\sum_{i,j=1}^{n}u^{ij}\frac{\partial}{\partial x_i}\wedge \frac{\partial}{\partial x_j}$ 
on a local coordinate $(x_1,\dots,x_n)$ of $M$. 
The Schouten bracket $[u, v ]_{Sch}$ for two 2-vectors $u$, $v$ is locally defined by 
\begin{eqnarray*}
[u, v ]_{Sch} &=&\sum_{i,j,k,l} [ u^{ij}\frac{\partial}{\partial x_i} \wedge\frac{\partial}{\partial x_j}, v^{kl}\frac{\partial}{\partial x_k} \wedge\frac{\partial}{\partial x_l} ]_{Sch}\\
&=&\sum_{i,j,k,l}\left\{ u^{ij}\frac{\partial v^{kl}}{\partial x_i}(\frac{\partial}{\partial x_j}\wedge\frac{\partial}{\partial x_k}\wedge\frac{\partial}{\partial x_l} ) - u^{ij}\frac{\partial v^{kl}}{\partial x_j}
(\frac{\partial}{\partial x_i}\wedge\frac{\partial}{\partial x_k}\wedge\frac{\partial}{\partial x_l} ) \right.\\
&&\left.+v^{kl}\frac{\partial  u^{ij}}{\partial x_k}(\frac{\partial}{\partial x_i}\wedge\frac{\partial}{\partial x_j}\wedge\frac{\partial}{\partial x_l} ) -v^{kl}\frac{\partial  u^{ij}}{\partial x_l}
(\frac{\partial}{\partial x_i}\wedge\frac{\partial}{\partial x_j}\wedge\frac{\partial}{\partial x_k}) \right\}
\end{eqnarray*}
where $u = \sum_{i,j=1}^{n}u^{ij}\frac{\partial}{\partial x_i}\wedge\frac{\partial}{\partial x_j}$ and 
$v=\sum_{k,l=1}^{n}v^{kl}\frac{\partial}{\partial x_k}\wedge \frac{\partial}{\partial x_l}$. 

We prepare a torsion free connection $\nabla$ of $M$. 
Then the Schouten bracket is also globally defined by the following map 
\[
\begin{array}{ccccc}
\Gamma(\wedge^2 TM \otimes \wedge^2 TM) & \longrightarrow & \Gamma(\wedge^2 TM \otimes TM)& \longrightarrow &\Gamma(\wedge^3 TM)  \\
\rotatebox{90}{$\in$}& &\rotatebox{90}{$\in$}& &\rotatebox{90}{$\in$} \\
(u, v) &\longmapsto & \nabla u\cdot v +\nabla v\cdot u &\longmapsto &\wedge^3(\nabla u\cdot v +\nabla v \cdot u)\\
\end{array}
\]
where $\nabla u\cdot v$ is the convolution of $\nabla u\in \Gamma (\wedge^2 TM \otimes T^* M)$ and $v\in \Gamma (\wedge^2 TM)$ 
and the last map $\wedge^3$ means the antisymmetrization. 

We consider a 2-vector $u$ as a homomorphism $T^*M \to TM$. 
Then the vector field $u(df)$ is the Hamiltonian vector field for a function $f$ on $M$. 
We can define a bracket $\{\cdot,\cdot \}$ by $\{f,g\}=u(df)g(=-u(dg)f)$ for functions $f,g$ on $M$. 
The bracket $\{\cdot,\cdot \}$ is a Poisson bracket on $M$, that is, satisfying the Jacobi identity, if $[u,u]_{Sch}=0$. 
We call that a 2-vector $u$ is a {\it Poisson structure} on $M$ if the Schouten bracket $[u, u ]_{Sch}$ vanishes. 
The equation $[u, u ]_{Sch}=0$ is written as a system of partial differential equations that is so hard to solve. 
%Globally, it is much harder to look for the solution of it. 
When the manifold is a complex manifold, the Schouten bracket is extended to complex 2-vectors by $\mathbb{C}$-linear. 
Then the Schouten bracket $[u, v ]_{Sch}$ for two holomorphic 2-vector $u$, $v$ is a holomorphic section of $\wedge^3 T^{1.0}M$. 
A holomorphic 2-vector $u$ is called a {\it holomorphic Poisson structure} if the Schouten bracket $[u, u ]_{Sch}$ vanishes. 
A holomorphic Poisson structure $u$ induces two (real) Poisson structures as its real part and imaginary part. 
On a complex two dimensional manifold, any holomorphic 2-vector is a holomorphic Poisson structure. 
However, it is even hard to find a holomorphic Poisson structure in the case that the complex dimension is greater than and equal to three.

\section{Holomorphic Poisson structures on $\mathbb{CP}^3$}
In this section, we study the space of holomorphic Poisson structures on $\mathbb{CP}^3$. 
We refer to~\cite{B} for Poisson structures on projective spaces. 

%Let $\widetilde{V}_q$ be the space of $\mathbb{C}^*$-invariant holomorphic $q$-vectors on $\mathbb{C}^{4}\backslash\{0\}$ :
%\[
%\widetilde{V}_q=\{ v \in \Gamma_{\text{$\mathbb{C}^*$-inv.}}(\wedge^qT^{\mathbb{C}}\mathbb{C}^{4}\backslash\{0\})\mid v: \text{holomorphic}\, \}.
%\]
%We define $V_q$ as the space of the push forward of elements of $\widetilde{V}_q$ by the projection $\pi: \mathbb{C}^{4}\backslash\{0\}\to \mathbb{CP}^{3}$ : 
%\[
%V_q=\{ d\pi(v) \in \Gamma(\wedge^qT^{\mathbb{C}}\mathbb{CP}^{3})\mid v\in \widetilde{V}_q \}.
%\]

\subsection{Holomorphic $q$-vectors on $\mathbb{CP}^3$}\label{s3.1}
Let $\pi$ be the projection $\mathbb{C}^4\backslash\{0\}\to \mathbb{CP}^3$. 
Any holomorphic $q$-vector on $\mathbb{CP}^3$ is given by the pushed forward 
of a $\mathbb{C}^*$-invariant holomorphic $q$-vector on $\mathbb{C}^4\backslash\{0\}$. 
We define $\widetilde{V}_q$ and $V_q$ as the space of $\mathbb{C}^*$-invariant holomorphic $q$-vectors on $\mathbb{C}^{4}\backslash\{0\}$ 
and that of holomorphic $q$-vectors on $\mathbb{CP}^3$, respectively. 
We denote by $(z_0, z_1, z_2, z_3)$ the standard coordinate of $\mathbb{C}^4$. 
Any holomorphic $1$-vector on $\mathbb{CP}^3$ is represented as the pushed forward of a holomorphic 1-vector $\sum_{i,j} a_{ij}z_i\frac{\partial}{\partial z_j}$ 
on $\mathbb{C}^4\backslash\{0\}$ since the coefficients are polynomials of degree 1. 
We consider the correspondence given by 
\begin{equation}\label{s3e2}
\sum_{i,j}a_{ij}z_i\frac{\partial}{\partial z_j} \mapsto (a_{ij}).
\end{equation}
Then we obtain the isomorphism  
\[
\widetilde{V}_1=\Bigl\{ \sum _{i,j}a_{ij}z_i\frac{\partial}{\partial z_j} \Bigm| a_{ij}\in \mathbb{C}  \Bigr\} \cong gl(4,\mathbb{C})
\]
where $gl(4,\mathbb{C})$ is the Lie algebra of the general linear group ${\rm GL}(4, \mathbb{C})$ of $\mathbb{C}^4$. 
The action of $\mathbb{C}^*$ on $\mathbb{C}^4\backslash\{0\}$ induces the vector field 
\[
\boldsymbol{l}=z_0\frac{\partial}{\partial z_0}+z_1\frac{\partial}{\partial z_1}+z_2\frac{\partial}{\partial z_2}+z_3\frac{\partial}{\partial z_3} 
\]
The 1-dimensional space $\left<\, \boldsymbol{l}\, \right>_{\mathbb{C}}$ is the space of holomorphic 1-vectors on $\mathbb{C}^4\backslash\{0\}$ 
which vanish by the push forward $d\pi$. 
Hence, the space $V_1$ of holomorphic 1-vectors on $\mathbb{CP}^3$ can be identified with 
the quotient space of $\widetilde{V}_1$ divided by $\left<\, \boldsymbol{l}\, \right>_{\mathbb{C}}$ : 
\[
V_1=\Bigl\{\sum _{i,j}a_{ij}z_i\frac{\partial}{\partial z_j} \Bigm| a_{ij}\in \mathbb{C} \Bigr\} \bigm/ \left<\, \boldsymbol{l}\, \right>_{\mathbb{C}}
\]
Under the correspondence (\ref{s3e2}), we obtain the isomorphism 
\[
V_1\cong gl(4,\mathbb{C})/\left< E \right>_{\mathbb{C}}\cong sl(4,\mathbb{C})
\]
where $sl(4,\mathbb{C})$ is the Lie algebra of the special linear group ${\rm SL}(4, \mathbb{C})$. 

%\subsection{Holomorphic 2-vectors on $\mathbb{CP}^3$}\label{s3.2}
A holomorphic 2-vector on $\mathbb{CP}^3$ is the pushed forward of a holomorphic 2-vector on $\mathbb{C}^4\backslash\{0\}$ whose coefficients are polynomials of degree 2. 
Such a 2-vector can be represented by $a_{ijkl}z_iz_j\frac{\partial}{\partial z_k}\wedge \frac{\partial}{\partial z_l}$ for a complex number $a_{ijkl}$. 
Then we can consider the coefficients $a_{ijkl}$ as the element $(a_{ijkl})$ of the tensor space $S^2\mathbb{C}^4\otimes \wedge^2(\mathbb{C}^4)^{\vee}$ 
where $S^2\mathbb{C}^4$ and $\wedge^2(\mathbb{C}^4)^{\vee}$ are the symmetric 2-tensor space of $\mathbb{C}^4$ and 
the anti-symmetric 2-tensor space of the dual space $(\mathbb{C}^4)^{\vee}$ of $\mathbb{C}^4$, respectively. 
Hence, 
\[
\widetilde{V}_2=\Bigl\{ a_{ijkl}z_iz_j\frac{\partial}{\partial z_k}\wedge \frac{\partial}{\partial z_l}  \Bigr\} \cong S^2\mathbb{C}^4\otimes \wedge^2(\mathbb{C}^4)^{\vee}.
\]
Since any holomorphic 2-vector which vanishes by $d\pi$ is represented by 
$(\sum_{i,j}a_{ij}z_i\frac{\partial}{\partial z_j})\wedge \boldsymbol{l}$ on $\mathbb{C}^4\backslash\{0\}$, %for $a_{ij}\in \mathbb{C}$
the space $V_2$ of holomorphic 2-vectors on $ \mathbb{CP}^3$ is regarded as the following quotient space : 
\[
V_2=\Bigl\{\sum_{i,j,k,l}a_{ijkl}z_iz_j\frac{\partial}{\partial z_k}\wedge \frac{\partial}{\partial z_l} \Bigm| a_{ijkl}\in \mathbb{C}\Bigr\} \Bigm/ 
\Bigl\{(\sum_{i,j}a_{ij}z_i\frac{\partial}{\partial z_j})\wedge \boldsymbol{l} \Bigm| a_{ij}\in \mathbb{C}\Bigr\}.
\]
%of the space $\{\sum_{i,j,k,l}a_{ijkl}z_iz_j\frac{\partial}{\partial z_k}\wedge \frac{\partial}{\partial z_l} \mid a_{ijkl}\in \mathbb{C}\}$ divided by 
%the subspace $\{(\sum_{i,j}a_{ij}z_i\frac{\partial}{\partial z_j})\wedge(\sum_{k}z_k\frac{\partial}{\partial z_k})\mid a_{ij}\in \mathbb{C}\}$. 
We also have the isomorphism 
\[
V_2\cong S^2\mathbb{C}^4\otimes \wedge^2(\mathbb{C}^4)^{\vee}\ / \ S^2\otimes \wedge^2(sl(4,\mathbb{C})\otimes \left< E \right>_{\mathbb{C}})
\]
where $S^2\otimes \wedge^2$ is the projection map 
from $\displaystyle gl(\mathbb{C}^4)\otimes gl(\mathbb{C}^4)\cong (\otimes^2\mathbb{C}^4)\otimes (\otimes^2 (\mathbb{C}^4)^{\vee})$ 
to $S^2\mathbb{C}^4\otimes \wedge^2(\mathbb{C}^4)^{\vee}$. 

As the same manner, the space $\widetilde{V}_3$ of $\mathbb{C}^*$-invariant $3$-vectors on $\mathbb{C}^{4}\backslash\{0\}$ is identified with 
\[
S^3\mathbb{C}^{4}\otimes \wedge^3(\mathbb{C}^{4})^{\vee}
\]
where $S^3\mathbb{C}^{4}$ and $\wedge^3(\mathbb{C}^{4})^{\vee}$ are the symmetric $3$-tensor space of $\mathbb{C}^{4}$ and 
the anti-symmetric $3$-tensor space of $(\mathbb{C}^{4})^{\vee}$, respectively. 
Then, the space $V_3$ of push forward of $\mathbb{C}^*$-invariant $3$-vectors on $\mathbb{C}^{4}\backslash\{0\}$ is isomorphic to 
\[
S^3\mathbb{C}^{4}\otimes \wedge^3(\mathbb{C}^{4})^{\vee}\, / \, S^3\otimes \wedge^3(\otimes^{2}sl(4,\mathbb{C})\otimes \left< E \right>_{\mathbb{C}})
\]
where $S^3\otimes \wedge^3$ is the projection 
from $\displaystyle \otimes^3 gl(\mathbb{C}^{4})\cong (\otimes^3\mathbb{C}^{4})\otimes (\otimes^3 (\mathbb{C}^{4})^{\vee})$ 
to $S^3\mathbb{C}^{4}\otimes \wedge^3(\mathbb{C}^{4})^{\vee}$. 
Hence we obtain the isomorphisms 
\begin{eqnarray*}
\widetilde{V}_q&\cong&S^q\mathbb{C}^{4}\otimes \wedge^q(\mathbb{C}^{4})^{\vee}, \label{s4.3e1} \\
V_q&\cong&S^q\mathbb{C}^{4}\otimes \wedge^q(\mathbb{C}^{4})^{\vee}\, / \, S^q\otimes \wedge^q(\otimes^{q-1}sl(4,\mathbb{C})\otimes \left< E \right>_{\mathbb{C}}) \notag
\end{eqnarray*}
for $q=1,2,3$. 
%by the argument in the previous section. 
%We will see the relation between $V_q$ and $V_q$ in the next subsection. 

\subsection{The Schouten bracket on $\mathbb{CP}^3$}\label{s3.3}
Let $(a_{ijkl})$ be an element of $S^2\mathbb{C}^{4}\otimes \wedge^2(\mathbb{C}^{4})^{\vee}$. 
By fixing $j$ and $l$, we can regard $(a_{ijkl})$ as the matrix whose $(i,k)$-component is $a_{ijkl}$, and denote the matrix by $A_{jl}$. 
Then, the matrix $(A_{jl})$ is the $gl(\mathbb{C}^{4})$-valued matrix with the $(j,l)$-component $A_{jl}$. 
Thus, $(a_{ijkl})$ is associated with the $gl(\mathbb{C}^{4})$-valued matrix $(A_{jl})$. 
If we take two elements $(a_{ijkl}), (a^{\prime}_{i^{\prime}j^{\prime}k^{\prime}l^{\prime}})$ of $S^2\mathbb{C}^{4}\otimes \wedge^2(\mathbb{C}^{4})^{\vee}$ 
and denote by $(A_{jl})$ and $(A^{\prime}_{j^{\prime}l^{\prime}})$ the associated matrices, respectively. 
Then we consider the bracket $[A_{jl}, A^{\prime}_{j^{\prime}l^{\prime}}]=A_{jl}A^{\prime}_{j^{\prime}l^{\prime}}-A^{\prime}_{j^{\prime}l^{\prime}}A_{jl}$ 
of matrices $A_{jl}, A^{\prime}_{j^{\prime}l^{\prime}}$ for each $j, l, j^{\prime}, l^{\prime}$ and 
have the $gl(\mathbb{C}^{4})$-valued matrix $([A_{jl}, A^{\prime}_{j^{\prime}l^{\prime}}])$ whose $(j,l)$-component is $[A_{jl}, A^{\prime}_{j^{\prime}l^{\prime}}]$. 
Moreover, $(([A_{jl}, A^{\prime}_{j^{\prime}l^{\prime}}]))$ is the $\otimes^2 gl(\mathbb{C}^{4})$-valued matrix whose $(j^{\prime},l^{\prime})$-component is $([A_{jl}, A^{\prime}_{j^{\prime}l^{\prime}}])$. 
We denote such $(([A_{jl}, A^{\prime}_{j^{\prime}l^{\prime}}]))$ by $[A_{jl}, A^{\prime}_{j^{\prime}l^{\prime}}]$ for simplicity. 
Thus, we can determine the element $[A_{jl}, A^{\prime}_{j^{\prime}l^{\prime}}]$ of $\otimes^3gl(\mathbb{C}^{4})$ 
for $(a_{ijkl}), (a^{\prime}_{i^{\prime}j^{\prime}k^{\prime}l^{\prime}})\in S^2\mathbb{C}^{4}\otimes \wedge^2(\mathbb{C}^{4})^{\vee}$ .

Under the identification of the space $\widetilde{V}_q$ of $\mathbb{C}^*$-invariant $q$-vectors on $\mathbb{C}^{4}\backslash\{0\}$ with 
$S^q\mathbb{C}^{4}\otimes \wedge^q(\mathbb{C}^{4})^{\vee}$, 
the Schouten bracket for $\mathbb{C}^*$-invariant 2-vectors on $\mathbb{C}^{4}\backslash\{0\}$ can be considered as the map 
\[
\widetilde{F} : S^2\mathbb{C}^{4}\otimes \wedge^2(\mathbb{C}^{4})^{\vee}\times S^2\mathbb{C}^{4}\otimes \wedge^2(\mathbb{C}^{4})^{\vee} \to S^3\mathbb{C}^{4}\otimes \wedge^3(\mathbb{C}^{4})^{\vee}
\]
defined by 
\[
\widetilde{F}(a_{ijkl}, a^{\prime}_{i^{\prime}j^{\prime}k^{\prime}l^{\prime}})=S^3\otimes \wedge^3(4[A_{jl}, A^{\prime}_{j^{\prime}l^{\prime}}])
\]
for $a_{ijkl}, a^{\prime}_{i^{\prime}j^{\prime}k^{\prime}l^{\prime}}\in S^2\mathbb{C}^{4}\otimes \wedge^2(\mathbb{C}^{4})^{\vee}$. 
Then, $\widetilde{F}$ is quadratic. 
We obtain 
\begin{equation}\label{s3e2.5}
[d\pi(v), d\pi(v^{\prime})]_{Sch}=d\pi[v, v^{\prime}]_{Sch}
\end{equation}
for $\mathbb{C}^*$-invariant 2-vectors $v, v^{\prime}$ on $\mathbb{C}^{4}\backslash\{0\}$. 
It follows from the $\pi$-relation (\ref{s3e2.5}) that the Schouten bracket on $\mathbb{CP}^3$ is given by 
the map from 
$S^2\mathbb{C}^{4}\otimes \wedge^2(\mathbb{C}^{4})^{\vee} / \{S^2\otimes \wedge^2(sl(4,\mathbb{C})\otimes \left< E \right>_{\mathbb{C}})\} 
\times S^2\mathbb{C}^{4}\otimes \wedge^2(\mathbb{C}^{4})^{\vee} / \{S^2\otimes \wedge^2(sl(4,\mathbb{C})\otimes \left< E \right>_{\mathbb{C}})\}$ 
to $S^3\mathbb{C}^{4}\otimes \wedge^3(\mathbb{C}^{4})^{\vee} / \{S^3\otimes \wedge^3(\otimes^{2}sl(4,\mathbb{C})\otimes \left< E \right>_{\mathbb{C}})\}$ induced by $\widetilde{F}$. 
We denote by $F$ the induced map. 
Then, we have 
\begin{prop}\label{s3.3p1}
The space of holomorphic Poisson structures on $\mathbb{CP}^3$ is identified with 
\[
\{a\in S^2\mathbb{C}^{4}\otimes \wedge^2(\mathbb{C}^{4})^{\vee}\, / \, S^2\otimes \wedge^2(sl(4,\mathbb{C})\otimes \left< E \right>_{\mathbb{C}}) \mid F(a,a)=0\}.
\]
$\hfill\Box$
\end{prop}

\subsection{The Schouten bracket on local coordinates}\label{s3.4}
We start to represent a holomorphic 2-vector $v= \sum_{i,j,k,l}a_{ijkl}z_iz_j\frac{\partial}{\partial z_k}\wedge \frac{\partial}{\partial z_l}$ on $\mathbb{C}^4\backslash\{0\}$ 
by using the local coordinate of $\mathbb{CP}^3$. 
Let $U_r$ be the open set $\{ [z_0,z_1,z_2,z_3]  \in  \mathbb{CP}^3 \mid z_r \ne 0 \}$ of $\mathbb{CP}^3$ for $r=0,1,2,3$. 
We denote by $\zeta_i$ the function $\zeta_i=\frac{z_i}{z_r}$ on $U_r$ for each $i=0, 1, 2, 3$. 
Then the system of functions $\{\zeta_0, \zeta_1, \zeta_2, \zeta_3\}\backslash \{\zeta_r\}$ without $\zeta_r$ defines the inhomogeneous coordinate on $U_r$. 
We remark that $\zeta_r=1$ on $U_r$. 
%Then $U_n$ corresponds to $\mathbb{C}^3$ by the correspondence $[ \frac{z_0}{z_n}  \cdot \cdot \cdot \frac{z_n}{z_n} \cdot \cdot\cdot  \frac{z_3}{z_n}]$ with $(t_1,t_2,t_3)$, that is $ \frac{z_i}{z_n} = t_{i'}$.
%For the 2-vector $v= \sum_{i,j,k,l}a_{ijkl}z_iz_j\frac{\partial}{\partial z_k}\wedge \frac{\partial}{\partial z_l}$ on $\mathbb{C}^4 \backslash \{ 0 \}$, 
The 2-vector $d\pi(v)$ is represented by 
\begin{equation}\label{s3e3}
d\pi(v)=\sum_{i,j,k\neq r,l\neq r}(a_{ijkl}\zeta_{i}\zeta_{j} - a_{ijkr}\zeta_{i}\zeta_{j}\zeta_{k} - a_{ijrl}\zeta_{i}\zeta_{j}\zeta_{l}) \frac{\partial }{\partial \zeta_{k}}\wedge \frac{\partial }{\partial \zeta_{l}}
\end{equation}
on the local coordinate $U_r$. 
For simplicity, we denote by $A_{r k l}$ the coefficient of $\frac{ \partial }{ \partial \zeta_{k}}\wedge  \frac { \partial }{ \partial \zeta_{l}}$ in (\ref{s3e3}) : 
\[
A_{r k l}=\sum_{i,j\neq r} a_{ijkl}\zeta_{i}\zeta_{j} - a_{ijkr}\zeta_{i}\zeta_{j}\zeta_{k} - a_{ijrl}\zeta_{i}\zeta_{j}\zeta_{l}.
\]
We represent the Schouten bracket of $d\pi(v)$ on the local coordinates of $\mathbb{CP}^3$. 
For each $r$, we can take a triple $(a,b,c)$ of $a,b,c,\in\{0,1,2,3\}\backslash \{r\}$ with $a<b<c$. 
Then the 3-vector $\frac{ \partial}{ \partial \zeta_a} \wedge \frac{ \partial}{ \partial \zeta_b} \wedge  \frac{ \partial}{\partial \zeta_c}$ 
generates the set of holomorphic 3-vectors on $U_r$, and it is denoted by $v_r$. 
Then the Schouten bracket $[d\pi(v),d\pi(v)]_{Sch}$ is written by 
\[
[d\pi(v),d\pi(v)]_{Sch}=2\sum_{(k,l,m)\in R_r}(\frac{ \partial A_{rkl}}{ \partial \zeta_l} A_{rlm} - A_{rkl} \frac{\partial A_{rlm}}{ \partial \zeta_l})v_r
\]
on each $U_r$, where $R_r$ be the set \{(a,b,c), (b,c,a), (c,a,b)\} of the rotation of the triple $(a,b,c)$. 
Hence the Schouten bracket $[d\pi(v),d\pi(v)]_{Sch}$ vanishes if and only if 
\[
\sum_{(k,l,m)\in R_r}(\frac{ \partial A_{rkl}}{ \partial \zeta_l}A_{rlm} - A_{rkl}\frac{ \partial A_{rlm}}{ \partial \zeta_l})=0
\]
on each $U_r$. It turns out that the space of holomorphic Poisson structures on $\mathbb{CP}^3$ is an algebraic variety. 
The space has the $\mathbb{C}^*$-action induced by the constant multiplication of 2-vectors. 
Then we obtain the following : 
\begin{prop}\label{s3.4p1}
The quotient of the space of holomorphic Poisson structures on $\mathbb{CP}^3$ by the $\mathbb{C}^*$-action 
is a projective algebraic variety. $\hfill\Box$
%\[
%\{a\in S^2\mathbb{C}^{4}\otimes \wedge^2(\mathbb{C}^{4})^{\vee}\, / \, S^2\otimes \wedge^2(sl(4,\mathbb{C})\otimes \left< E \right>_{\mathbb{C}}) \mid F(a,a)=0\}\, / \,\mathbb{C}^*.
%\]
\end{prop}

\section{Poisson structures on $S^4$}
Let $\mathbb{H}$ be the Hamilton's quaternionic number field $\mathbb{R}+i\mathbb{R}+j\mathbb{R}+k\mathbb{R}$ . 
Then $\mathbb{H}^{2}$ is isomorphic to the right $\mathbb{C}$-module $\mathbb{C}^{4}$ under the identification of $\mathbb{H}$ with $\mathbb{C}^2 (\cong \mathbb{C} + j \mathbb{C})$. 
Let $p$ be the projection from $\mathbb{C}^{4} \backslash\{0\}\cong \mathbb{H}^{2}\backslash\{0\}$ to $\mathbb{HP}^1=S^4$. 
In this section, we consider real $q$-vectors on $\mathbb{HP}^1$ given by the push forward of $\mathbb{H}^*$-invariant real $q$-vectors with holomorphic $(q,0)$-parts on $\mathbb{C}^{4}\backslash\{0\}$. 

\subsection{$\mathbb{H}^*$-invariant real $q$-vectors on $\mathbb{C}^{4}\backslash\{0\}$}\label{s4.1}
Let $\widetilde{W}_q$ be the space of $\mathbb{H}^*$-invariant real $q$-vectors with holomorphic $(q,0)$-parts on $\mathbb{C}^{4}\backslash\{0\}$ :
\[
\widetilde{W}_q=\{ w \in \Gamma_{\text{$\mathbb{H}^*$-inv.}}(\wedge^qT\mathbb{C}^{4}\backslash\{0\})\mid w^{q,0}: \text{holomorphic}\, \}.
\] 
The right action of $j$ on $\mathbb{C}^{4} \backslash\{0\}$ induces the map $\boldsymbol{j} : \mathbb{C}^{4}\backslash\{0\} \to \mathbb{C}^{4}\backslash \{0\}$ given by 
\[
(z_0,z_1,z_2,z_3) \mapsto (- \bar{z}_1,\bar{z}_0,-\bar{z}_3,\bar{z}_2)
\]
for $(z_0,z_1,z_2,z_3)\in \mathbb{C}^{4}\backslash\{0\}$. 
Let $\boldsymbol{l}'$ denote the complex vector field on $\mathbb{C}^{4} \backslash\{0\}$ generated by $\exp(tj)$ for $t\in \mathbb{R}$. 
Then we have 
\[
\boldsymbol{l}'=\bar{z}_0\frac{\partial}{\partial z_1}-\bar{z}_1\frac{\partial}{\partial z_0}+\bar{z}_{2}\frac{\partial}{\partial z_{3}}-\bar{z}_{3}\frac{\partial}{\partial z_{2}}.
\]
Amy $\mathbb{C}^*$-invariant real $q$-vector $w$ on $\mathbb{C}^{4}\backslash\{0\}$ is $\mathbb{H}^*$-invariant 
if and only if it is preserved by $\boldsymbol{l}'$, that is, $L_{(\boldsymbol{l}'+\overline{\boldsymbol{l}'})}w=L_{(\boldsymbol{l}'-\overline{\boldsymbol{l}'})}w=0$. 
The condition $L_{(\boldsymbol{l}'+\overline{\boldsymbol{l}'})}w=L_{(\boldsymbol{l}'-\overline{\boldsymbol{l}'})}w=0$ is equal to $L_{\boldsymbol{l}'}w=0$ since $w$ is real. 
Hence, the space $\widetilde{W}_q$ is written by 
\[
\widetilde{W}_q=\{ w \in \Gamma_{\text{$\mathbb{C}^*$-inv.}}(\wedge^qT\mathbb{C}^{4}\backslash\{0\})\mid L_{\boldsymbol{l}'}w=0,\ w^{q,0}: \text{holomorphic}\, \}. 
\]
%where $L_{\boldsymbol{l}'}w$ means the Lie derivative of $w$ by $\boldsymbol{l}'$. 
We also define $W_q$ by the space of real $q$-vectors on $\mathbb{HP}^1$ induced by elements of $\widetilde{W}_q$ : 
\[
W_q=\{ dp(w) \in \Gamma(\wedge^qT\mathbb{HP}^1)\mid w\in \widetilde{W}_q \}.
\] 
The real and imaginary parts of $\boldsymbol{l}$ and $\boldsymbol{l}'$ span the real 4-dimensional vector space. 
Then the intersection of the space and $\widetilde{W}_1$ is the 1-dimensional vector space generated by $\boldsymbol{l}+\bar{\boldsymbol{l}}$. 
Thus the space $\left< \boldsymbol{l}+\bar{\boldsymbol{l}} \right>_{\mathbb{R}}$ is the kernel of the push forward by $p$, 
and $W_1$ is identified with the the quotient space $\widetilde{W}_1/\left< \boldsymbol{l}+\bar{\boldsymbol{l}} \right>_{\mathbb{R}}$. 
Moreover, there exists an isomorphism 
\[
W_q\cong \widetilde{W}_q\, /\, \widetilde{W}_{q-1}\wedge \left< \boldsymbol{l}+\bar{\boldsymbol{l}} \right>_{\mathbb{R}}
\]
where the right hand side is the quotient space of $\widetilde{W}_q$ divided by $\widetilde{W}_{q-1}\wedge \left< \boldsymbol{l}+\bar{\boldsymbol{l}} \right>_{\mathbb{R}}$. 
%where $\boldsymbol{l}$ is the vector field $\sum_{i=0}^{3}z_i\frac{\partial}{\partial z_i}$ on $\mathbb{C}^{4}\backslash\{0\}$. 

\subsection{$q$-vectors on $S^4$}\label{s4.2}
We consider the space $\widetilde{W}_1$ of $\mathbb{H}^*$-invariant real $1$-vectors on $\mathbb{C}^{4}\backslash\{0\}$, whose $(1,0)$-parts are holomorphic. 
%We consider a 1-vector on $\mathbb{HP}^1$ corresponding to a $\mathbb{C}^*$-invariant 1-vector on $\mathbb{C}^{4}\backslash\{0\}$ preserved by $\boldsymbol{l}'$. 
The space of $\mathbb{C}^*$-invariant real 1-vectors on $\mathbb{C}^{4} \backslash\{0\}$ with holomorphic $(1,0)$-parts is 
\[
\Bigl\{ \sum_{i,j} (a_{ij}z_i\frac{\partial}{\partial z_j} + \bar{a}_{ij}\bar{z}_i\frac{\partial}{\partial {\bar{z}_j} }) \Bigm| a_{ij}\in \mathbb{C}\Bigr\}.
\]
%which is the space of the real part of a holomorphic 1-vector $\sum_{i,j} a_{ij}z_i\frac{\partial}{\partial z_j}$ on $\mathbb{C}^{4} \backslash\{0\}$. 
Let $\varphi$ be the mapping $(0,1,2, 3) \to (1,0,3, 2)$, that is, $\varphi(0)=1, \varphi(1)=0, \varphi(2)=3, \varphi(3)=2$. 
Now, 
\[
L_{\boldsymbol{l}'}z_i=(-1)^{\varphi(i)}\bar{z}_{\varphi(i)},\quad L_{\boldsymbol{l}'}\bar{z}_i=0,\quad 
L_{\boldsymbol{l}'}\frac{\partial}{\partial z_j}=0,\quad  L_{\boldsymbol{l}'}\frac{\partial}{\partial \bar{z}_j}=(-1)^{\varphi(j)}\frac{\partial}{\partial z_{\varphi(j)}}
\]
for each $i,j$. 
It yields that  
\begin{eqnarray*}
L_{\boldsymbol{l}'}\left( \sum_{i,j} a_{ij}z_i\frac{\partial}{\partial z_j} + \bar{a}_{ij}\bar{z}_i\frac{\partial}{\partial {\bar{z}_j} } \right)
&=&\sum_{i,j} (-1)^{\varphi(i)}a_{ij}\bar{z}_{\varphi(i)}\frac{\partial}{\partial z_j} + (-1)^{\varphi(j)}\bar{a}_{ij}\bar{z}_i\frac{\partial}{\partial {\bar{z}_{\varphi(j)}}} \\
&=&\sum_{i,j} ((-1)^{\varphi(i)}a_{ij} + (-1)^{j}\bar{a}_{\varphi(i)\varphi(j)})\bar{z}_{\varphi(i)}\frac{\partial}{\partial z_j}. 
\end{eqnarray*}
Hence, the 1-vector $\sum_{i,j} (a_{ij}z_i\frac{\partial}{\partial z_j}+ \bar{a}_{ij}\bar{z}_i\frac{\partial}{\partial {\bar{z}_j}})$ is $\mathbb{H}^*$-invariant if and only if 
\begin{equation}\label{s4e1}
a_{ij}=(-1)^{i+j}\bar{a}_{\varphi(i)\varphi(j)}.
\end{equation}
Thus the space $\widetilde{W}_1$ is written by 
\[
\widetilde{W}_1=\Bigl\{ \sum_{i,j} (a_{ij}z_i\frac{\partial}{\partial z_j}+ \bar{a}_{ij}\bar{z}_i\frac{\partial}{\partial {\bar{z}_j}}) 
\Bigm| a_{ij}=(-1)^{i+j}\bar{a}_{\varphi(i)\varphi(j)}\Bigr\}.
\]
We can consider the coefficients $a_{ij}$ as the matrix $(a_{ij})$ of $gl(4,\mathbb{C})$ under the correspondence 
\begin{equation*}\label{s4e2}
\sum_{i,j}(a_{ij}z_i\frac{\partial}{\partial z_j}+ \bar{a}_{ij}\bar{z}_i\frac{\partial}{\partial {\bar{z}_j}}) \mapsto (a_{ij}).
\end{equation*}
Let ${\rm GL}(2, \mathbb{H})$ be the general linear group of $\mathbb{H}^{2}$. 
Then the Lie algebra $gl(2, \mathbb{H})$ of ${\rm GL}(2, \mathbb{H})$ is identified with the subspace of $gl(4, \mathbb{C})$ by the map 
\begin{equation}\label{s4e2.5}
\begin{pmatrix}
a_{00}+jb_{00}&a_{01}+jb_{01}\\
a_{10}+jb_{10}&a_{11}+jb_{11}
\end{pmatrix}
\longmapsto 
\begin{pmatrix}
\left.
\begin{array}{cc}
a_{00} & -\bar{b}_{00}\\
b_{00} & \bar{a}_{00}
\end{array}
\right.
\left.
\begin{array}{cc}
a_{01} & -\bar{b}_{01}\\
b_{01} & \bar{a}_{01}
\end{array}
\right. \\
\left.
\begin{array}{cc}
a_{10} & -\bar{b}_{10}\\
b_{10} & \bar{a}_{10}
\end{array}
\right.
\left.
\begin{array}{cc}
a_{11} & -\bar{b}_{11}\\
b_{11} & \bar{a}_{11}
\end{array}
\right.
\end{pmatrix}
\end{equation}
for $(a_{kl}+jb_{kl})\in gl(2, \mathbb{H})$. 
Conversely, an element $(a_{ij})$ of $gl(4,\mathbb{C})$ is in the image of $gl(2, \mathbb{H})$ by the map (\ref{s4e2.5}) if $a_{ij}$ satisfies the equation (\ref{s4e1}). 
%The condition (\ref{s4e1}) is equivalent that the matrix $(a_{ij})\in gl(4,\mathbb{C})$ is contained in 
%the Lie algebra $gl(2, \mathbb{H})$ of the general linear group ${\rm GL}(2, \mathbb{H})$ of $\mathbb{H}^2$. 
%under the identification ${\rm GL}(2, \mathbb{H})$ 
%We denote by $W$ the space of $\mathbb{H}^*$-invariant real 1-vectors on $\mathbb{C}^{4} \backslash\{0\}$. 
Hence, we obtain 
\[
\widetilde{W}_1\cong gl(2, \mathbb{H}).
\]
%Secondly, we consider the space $W_1$ of real $1$-vectors on $\mathbb{HP}^1$ induced by elements of $\widetilde{W}_1$. 
%and define $W$ the space of such elements of $W$. 
%The real and imaginary parts of $\boldsymbol{l}$ and $\boldsymbol{l}'$ span the real 4-dimensional vector space. 
%Then the intersection of the space and $\widetilde{W}_1$ is the 1-dimensional vector space generated by $\boldsymbol{l}+\bar{\boldsymbol{l}}$. 
%Thus the space $\left< \boldsymbol{l}+\bar{\boldsymbol{l}} \right>_{\mathbb{R}}$ is the kernel of the push forward by $p$, 
%and $W_1$ is identified with the the quotient space $\widetilde{W}_1/\left< \boldsymbol{l}+\bar{\boldsymbol{l}} \right>_{\mathbb{R}}$. 
Since the real 1-vector $\boldsymbol{l}+\bar{\boldsymbol{l}}$ corresponds to the unit matrix $E \in gl(2, \mathbb{H})$, %by the map (\ref{s4e2}), 
we also have 
\[
W_1\cong gl(2, \mathbb{H})/\left< E \right>_{\mathbb{R}}. %\cong sl(2, \mathbb{H})
\]
%of $gl(2, \mathbb{H})$ divided by the 1-dimensional vector space $\left< E \right>_{\mathbb{R}}$. 
%where $sl(2, \mathbb{H})$ is the Lie algebra of the special linear group ${\rm SL}(2, \mathbb{H})$. 

%\subsection{2-vectors on $\mathbb{HP}^1$}\label{s4.3}
We consider the space $\widetilde{W}_2$ of $\mathbb{H}^*$-invariant real $2$-vectors on $\mathbb{C}^{4}\backslash\{0\}$ with holomorphic $(2,0)$-parts. 
We take a $\mathbb{C}^*$-invariant real 2-vector $w$ on $\mathbb{C}^{4} \backslash\{0\}$ given by 
\[
w=\sum_{i,j,k,l}a_{ijkl}z_iz_j\frac{\partial}{\partial z_k}\wedge \frac{\partial}{\partial z_l}
+b_{ijkl}z_i\bar{z}_j\frac{\partial}{\partial z_k}\wedge \frac{\partial}{\partial \bar{z}_l}
+\bar{a}_{ijkl}\bar{z}_i\bar{z}_j\frac{\partial}{\partial \bar{z}_k}\wedge \frac{\partial}{\partial \bar{z}_l}
\]
for $a_{ijkl}, b_{ijkl}\in \mathbb{C}$. Then the equations $a_{jikl}=a_{ijkl},\ a_{ijlk}=-a_{ijkl},\ b_{jilk}=-\bar{b}_{ijkl}$ hold for $i,j,k,l=0,\dots,3$. 
We recall that the coefficients $a_{ijkl}$ define the element $(a_{ijkl})$ of $S^2\mathbb{C}^{4}\otimes \wedge^2(\mathbb{C}^{4})^{\vee}$ as in \S\ref{s3.1}. 
Now, 
\begin{eqnarray*}
L_{\boldsymbol{l}'}w
&=&\sum_{i,j,k,l} a_{ijkl}\left( (-1)^{\varphi(i)}\bar{z}_{\varphi(i)}z_j + (-1)^{\varphi(j)}z_i \bar{z}_{\varphi(j)}\right) \frac{\partial}{\partial z_k}\wedge \frac{\partial}{\partial z_l} \\
&&+ b_{ijkl}\left( (-1)^{\varphi(i)}\bar{z}_{\varphi(i)}\bar{z}_j \frac{\partial}{\partial z_k}\wedge \frac{\partial}{\partial \bar{z}_l} 
+(-1)^{\varphi(l)}z_i \bar{z}_{j}\frac{\partial}{\partial z_k}\wedge \frac{\partial}{\partial z_{\varphi(l)}}\right) \\
&&+ \bar{a}_{ijkl}\bar{z}_{i}\bar{z}_j\left( (-1)^{\varphi(k)} \frac{\partial}{\partial z_{\varphi(k)}}\wedge \frac{\partial}{\partial \bar{z}_l} 
+(-1)^{\varphi(l)}\bar{z}_i \bar{z}_{j}\frac{\partial}{\partial \bar{z}_k}\wedge \frac{\partial}{\partial z_{\varphi(l)}}\right) \\
%&=&\sum_{i,j,k,l} 2a_{ijkl}(-1)^{\varphi(j)}z_i \bar{z}_{\varphi(j)}\frac{\partial}{\partial z_k}\wedge \frac{\partial}{\partial z_l} \\
%&&+ b_{ijkl}\left( (-1)^{\varphi(i)}\bar{z}_{\varphi(i)}\bar{z}_j \frac{\partial}{\partial z_k}\wedge \frac{\partial}{\partial \bar{z}_l} 
%+(-1)^{\varphi(l)}z_i \bar{z}_{j}\frac{\partial}{\partial z_k}\wedge \frac{\partial}{\partial z_{\varphi(l)}}\right) \\
%&&+ 2\bar{a}_{ijkl}(-1)^{\varphi(k)}\bar{z}_{i}\bar{z}_j \frac{\partial}{\partial z_{\varphi(k)}}\wedge \frac{\partial}{\partial \bar{z}_l} \\
%&=&\sum_{i,j,k,l} \left( 2a_{ijk\varphi(l)}(-1)^{\varphi(j)}z_i \bar{z}_{\varphi(j)}+b_{ijkl}(-1)^{\varphi(l)}z_i \bar{z}_{j}\right) \frac{\partial}{\partial z_k}\wedge \frac{\partial}{\partial z_{\varphi(l)}} \\
%&&+ \left( b_{ijkl}(-1)^{\varphi(i)}\bar{z}_{\varphi(i)}\bar{z}_j+ 2\bar{a}_{ij\varphi(k)l}(-1)^{k}\bar{z}_{i}\bar{z}_j \right) 
%\frac{\partial}{\partial z_k}\wedge \frac{\partial}{\partial \bar{z}_l} \\
&=&\sum_{i,j,k,l} \left( 2a_{i\varphi(j)k\varphi(l)}(-1)^{j}+b_{ijkl}(-1)^{\varphi(l)}\right) z_i \bar{z}_{j} \frac{\partial}{\partial z_k}\wedge \frac{\partial}{\partial z_{\varphi(l)}} \\
&&+ \left( b_{ijkl}(-1)^{\varphi(i)}+ 2\bar{a}_{\varphi(i)j\varphi(k)l}(-1)^{k} \right) 
\bar{z}_{\varphi(i)}\bar{z}_j \frac{\partial}{\partial z_k}\wedge \frac{\partial}{\partial \bar{z}_l}. \\
\end{eqnarray*}
Hence, the $\mathbb{H}^*$-invariant condition of $w$ is equal to 
\begin{equation*}\label{s4e3}
b_{ijkl}=2(-1)^{j+l}a_{i\varphi(j)k\varphi(l)}=2(-1)^{i+k}\bar{a}_{\varphi(i)j\varphi(k)l}.
\end{equation*}
It yields that $b_{ijkl}$ is uniquely determined by $(a_{ijkl})$ satisfying 
\begin{equation}\label{s4e4}
a_{ijkl}=(-1)^{i+j+k+l}\bar{a}_{\varphi(i)\varphi(j)\varphi(k)\varphi(l)}. 
\end{equation}
The space of elements of $S^2\mathbb{C}^{4}\otimes \wedge^2(\mathbb{C}^{4})^{\vee}$ satisfying the condition (\ref{s4e4}) 
is the tensor space $S^2\otimes \wedge^2(gl(2, \mathbb{H})\otimes gl(2, \mathbb{H}))$. 
%\[
%W=\{(a_{ijkl})\in S^2\mathbb{C}^{4}\otimes \wedge^2(\mathbb{C}^{4})^{\vee} \mid a_{ijkl}=(-1)^{i+j+k+l}\bar{a}_{\varphi(i)\varphi(j)\varphi(k)\varphi(l)}\}.
%\]
%Then the space $\widetilde{W}_2$ is isomorphic to $V$ 
Hence, by the correspondence %map from the space of 2-vectors on $\mathbb{C}^{4} \backslash\{0\}$ to $S^2\mathbb{C}^{4}\otimes \wedge^2(\mathbb{C}^{4})^{\vee}$ given by
\begin{equation*}\label{s4e5}
\sum_{i,j,k,l}a_{ijkl}z_iz_j\frac{\partial}{\partial z_k}\wedge \frac{\partial}{\partial z_l}
+b_{ijkl}z_i\bar{z}_j\frac{\partial}{\partial z_k}\wedge \frac{\partial}{\partial \bar{z}_l}
+\bar{a}_{ijkl}\bar{z}_i\bar{z}_j\frac{\partial}{\partial \bar{z}_k}\wedge \frac{\partial}{\partial \bar{z}_l} \mapsto (a_{ijkl}), 
\end{equation*}
we obtain the isomorphism 
\[
\widetilde{W}_2\cong S^2\otimes \wedge^2(gl(2, \mathbb{H})\otimes gl(2, \mathbb{H})).
\]
It yields that %the space $W_2$ of real 2-vectors on $\mathbb{HP}^1$ induced by $\widetilde{W}_2$ is isomorphic the quotient space 
\[
W_2\cong S^2\otimes \wedge^2(gl(2, \mathbb{H})\otimes gl(2, \mathbb{H}))\, / \, S^2\otimes \wedge^2(gl(2, \mathbb{H})\otimes \left< E \right>_{\mathbb{R}}).
\]
%of $V$ divided by the vector space $S^2\otimes \wedge^2(gl(2, \mathbb{H})\otimes \left< E \right>_{\mathbb{R}})$. 

Let $w$ be a $\mathbb{C}^*$-invariant real 3-vector on $\mathbb{C}^{4} \backslash\{0\}$. 
Then $w$ is given by 
\begin{eqnarray*}
w&=&\sum_{i,j,k,l,m,n}a_{ijklmn}z_iz_jz_k\frac{\partial}{\partial z_l}\wedge \frac{\partial}{\partial z_m}\wedge \frac{\partial}{\partial z_n}
+b_{ijklmn}z_iz_j\bar{z}_k\frac{\partial}{\partial z_l}\wedge \frac{\partial}{\partial z_m}\wedge \frac{\partial}{\partial \bar{z}_n} \\
&&+\bar{b}_{ijklmn}\bar{z}_i\bar{z}_jz_k\frac{\partial}{\partial \bar{z}_l}\wedge \frac{\partial}{\partial \bar{z}_m}\wedge \frac{\partial}{\partial z_n}
+\bar{a}_{ijklmn}\bar{z}_i\bar{z}_j\bar{z}_k\frac{\partial}{\partial \bar{z}_l}\wedge \frac{\partial}{\partial \bar{z}_m}\wedge \frac{\partial}{\partial \bar{z}_n}
\end{eqnarray*}
for $a_{ijklmn}, b_{ijklmn}\in \mathbb{C}$. 
Then $a_{\sigma(j)\sigma(i)\sigma(k)lmn}=a_{ijklmn},\ a_{ijk\tau(l)\tau(m)\tau(n)}={\rm sgn}(\tau) a_{ijklmn}$ for permutations $\sigma$ of $\{i,j,k\}$ and $\tau$ of $\{l,m,n\}$ and 
$b_{jiklmn}=\bar{b}_{ijklmn},\ b_{ijkmln}=-\bar{b}_{ijklmn}$ for $i,j,k,l,m,n=0,\dots,3$. 
We can find that the space $\widetilde{W}_3$ is identified 
with the tensor space $S^3\otimes \wedge^3(\otimes^3 gl(2, \mathbb{H}))$ : 
\[
\widetilde{W}_3\cong S^3\otimes \wedge^3(\otimes^3 gl(2, \mathbb{H})).
\]
We have 
\begin{eqnarray*}
L_{\boldsymbol{l}'}w
&=&\sum_{i,j,k,l,m,n} 3a_{ijklmn}(-1)^{\varphi(k)}z_i z_j \bar{z}_{\varphi(k)}\frac{\partial}{\partial z_l}\wedge \frac{\partial}{\partial z_m}\wedge \frac{\partial}{\partial z_n} \\
&&+b_{ijklmn}\left( 2(-1)^{\varphi(j)}z_i\bar{z}_{\varphi(j)}\bar{z}_k \frac{\partial}{\partial z_l}\wedge \frac{\partial}{\partial z_m}\wedge \frac{\partial}{\partial \bar{z}_n} 
+(-1)^{\varphi(n)}z_i z_j \bar{z}_{k}\frac{\partial}{\partial z_l}\wedge \frac{\partial}{\partial z_m}\wedge \frac{\partial}{\partial z_{\varphi(n)}}\right) \\
&&+\bar{b}_{ijklmn}\left( (-1)^{\varphi(k)}\bar{z}_i\bar{z}_j\bar{z}_{\varphi(k)} \frac{\partial}{\partial z_l}\wedge \frac{\partial}{\partial z_m}\wedge \frac{\partial}{\partial z_n} 
+(-1)^{\varphi(m)}\bar{z}_i \bar{z}_j z_k\frac{\partial}{\partial \bar{z}_l}\wedge \frac{\partial}{\partial z_{\varphi(m)}}\wedge \frac{\partial}{\partial z_n}\right) \\
&&+3\bar{a}_{ijklmn}(-1)^{\varphi(n)}\bar{z}_i \bar{z}_j \bar{z}_k \frac{\partial}{\partial \bar{z}_l}\wedge \frac{\partial}{\partial \bar{z}_m}\wedge \frac{\partial}{\partial z_{\varphi(n)}} \\
&=&\sum_{i,j,k,l,m,n} \left( 3a_{ij\varphi(k)lm\varphi(n)}(-1)^{k}+b_{ijklmn}(-1)^{\varphi(n)}\right) 
z_i z_j \bar{z}_{k}\frac{\partial}{\partial z_l}\wedge\frac{\partial}{\partial z_m}\wedge \frac{\partial}{\partial z_{\varphi(n)}} \\
&&+2\left( -b_{k\varphi(j)in\varphi(m)l}(-1)^{j}+ \bar{b}_{ijklmn}(-1)^{\varphi(m)}\right) 
\bar{z}_i\bar{z}_j z_k \frac{\partial}{\partial \bar{z}_l} \wedge\frac{\partial}{\partial z_{\varphi(m)}}\wedge\frac{\partial}{\partial z_n}\\
&&+\left( 3\bar{a}_{ij\varphi(k)lm\varphi(n)}(-1)^{n}+ \bar{b}_{ijklmn}(-1)^{\varphi(k)} \right) 
\bar{z}_i\bar{z}_j\bar{z}_{\varphi(k)} \frac{\partial}{\partial \bar{z}_l}\wedge \frac{\partial}{\partial \bar{z}_m}\wedge \frac{\partial}{\partial z_n}. 
\end{eqnarray*}
It implies that the $\mathbb{H}^*$-invariant condition of $w$ is given by 
\begin{equation*}\label{s4e6}
b_{ijklnm}=(-1)^{j+m}\bar{b}_{i\varphi(j)kl\varphi(m)n}=(-1)^{k+n}3a_{ij\varphi(k)lm\varphi(n)}.
\end{equation*}
Hence, $b_{ijklmn}$ is uniquely determined by $(a_{ijklmn})$ satisfying 
\begin{equation}\label{s4e7}
a_{ijklmn}=(-1)^{i+j+k+l+m+n}\bar{a}_{\varphi(i)\varphi(j)\varphi(k)\varphi(l)\varphi(m)\varphi(n)}. 
\end{equation}
We also have the isomorphism  
\[
W_3\cong S^3\otimes \wedge^3(\otimes^3 gl(2, \mathbb{H}))\, / \, S^3\otimes \wedge^3(\otimes^{2}gl(2,\mathbb{H})\otimes \left< E \right>_{\mathbb{R}}).
\]
In the next subsection, we will see that the space $W_q$ is identified with a real form of the space induced by 
$\mathbb{C}^*$-invariant holomorphic $q$-vectors on $\mathbb{C}^{4}\backslash\{0\}$ for each $q=1,2,3$.

\subsection{The real structures on $\widetilde{V}_q$ and $V_q$}\label{s4.4}
We define the map $\Phi: \Gamma(T^{\mathbb{C}}\mathbb{C}^{4}\backslash\{0\}) \to \Gamma(T^{\mathbb{C}}\mathbb{C}^{4}\backslash\{0\})$ by 
\[
\Phi(u)=\overline{d\boldsymbol{j}(u)}
\]
for any $u\in \Gamma(T^{\mathbb{C}}\mathbb{C}^{4}\backslash\{0\})$. 
Then the map $\Phi$ is anti-$\mathbb{C}$-linear and preserves the inner product given by the standard Riemannian metric on $\mathbb{C}^{4}$. 
Since $\Phi^2(u)=u$ for any $\mathbb{C}^{*}$-invariant holomorphic 1-vector $u$ on $\mathbb{C}^{4} \backslash\{0\}$, 
the map $\Phi$ defines the real structure on the space $\widetilde{V}_q$ of $\mathbb{C}^{*}$-invariant holomorphic $q$-vectors on $\mathbb{C}^{4} \backslash\{0\}$ for each $q$. 
We denote by $\widetilde{V}_q^{\rm Re}$ the real form of $\Phi$, that is, the set of fixed points of $\Phi$. 
%Then the space $\widetilde{W}_1$ of $\mathbb{H}^*$-invariant real $1$-vectors on $\mathbb{C}^{4}\backslash\{0\}$ 
%with holomorphic $(1,0)$-parts is isomorphic to $\widetilde{V}_1^{\rm Re}$ by the map $w\mapsto w^{1,0}$ for $w\in \widetilde{W}_1$. 
%The map $\Phi$ induces the real structure on the space $\widetilde{V}_q$ of $\mathbb{C}^{*}$-invariant holomorphic $q$-vectors on $\mathbb{C}^{4} \backslash\{0\}$, 
%which is denoted by $\Phi$ for simplicity. 
%Then the space $\widetilde{W}_q$ is isomorphic to the real form $\widetilde{V}_q^{\rm Re}$ of $\Phi$ 
\begin{prop}\label{s4.4p1}
There exists an isomorphism $\widetilde{W}_q\cong \widetilde{V}_q^{\rm Re}$ given by the map 
\begin{equation*}\label{s4.4e1}
w\mapsto w^{q,0}
\end{equation*}
for $w\in \widetilde{W}_q$. 
\end{prop}
\begin{proof}
Let $w$ be an element of $\widetilde{W}_q$. 
Then $w$ is $\mathbb{H}^*$-invariant if and only if the coefficient matrix $(a_{i_1 i_2 \cdots i_{2q-1} i_{2q}})$ of the $(q,0)$-part $w^{q,0}$ satisfies 
\[
a_{i_1 i_2 \cdots i_{2q}}=(-1)^{i_1+i_2+\dots+i_{2q}}\bar{a}_{\varphi(i_1) \varphi(i_2) \cdots \varphi(i_{2q})}
\]
by the equations (\ref{s4e1}), (\ref{s4e4}) and (\ref{s4e7}) in the previous subsection. 
On the other hand, under the identification of $\widetilde{V}_q\cong S^q\mathbb{C}^{4}\otimes \wedge^q(\mathbb{C}^{4})^{\vee}$, 
the real structure $\Phi$ can be written by 
\[
\Phi(a_{i_1 i_2 \cdots i_{2q}})=(-1)^{i_1+i_2+\dots+i_{2q}}\bar{a}_{\varphi(i_1) \varphi(i_2) \cdots \varphi(i_{2q})}
\]
for $(a_{i_1 i_2 \cdots i_{2q-1} i_{2q}})\in S^q\mathbb{C}^{4}\otimes \wedge^q(\mathbb{C}^{4})^{\vee}$. 
Hence, $w^{q,0}$ is in $\widetilde{V}_q^{\rm Re}$. 
It means that the correspondence $w\mapsto w^{q,0}$ defines the map $\widetilde{W}_q\to \widetilde{V}_q^{\rm Re}$. 
Conversely, for any element $w^{q,0}$ of $\widetilde{V}_q^{\rm Re}$, 
we can uniquely determined an element $w$ of $\widetilde{W}_q$ such that the $(q,0)$-part is $w^{q,0}$ by the argument in  previous subsection. 
Hence the map $\widetilde{W}_q\to \widetilde{V}_q^{\rm Re}$ is isomorphic. 
\end{proof}
It follows from $\Phi( \boldsymbol{l} )= \boldsymbol{l}$ that $\Phi$ preserves the complex vector space $\left< \boldsymbol{l} \right>_{\mathbb{C}}$. 
Hence, $\Phi$ provides the real structure on the space 
$V_q \cong \widetilde{V}_q\, /\, \widetilde{V}_{q-1}\wedge \left< \boldsymbol{l} \right>_{\mathbb{C}}$. 
We denote by $\Phi$ the real structure on $V_q$, again. 
Then, we regard the space $W_q$ as the real form $V_q^{\rm Re}$ of $V_q$ :
\begin{cor}\label{s4.4c1}
There exists an isomorphism $W_q\cong V_q^{\rm Re}$ given by the map 
\begin{equation*}
dp(w)\mapsto d\pi(w^{q,0})
\end{equation*}
for $w\in \widetilde{W}_q$. 
\end{cor}
\begin{proof} 
$W_q\cong \widetilde{W}_q\, /\, \widetilde{W}_{q-1}\wedge \left< \boldsymbol{l}+\bar{\boldsymbol{l}} \right>_{\mathbb{R}} 
\cong \widetilde{V}_q^{\rm Re}\, /\, \widetilde{V}_{q-1}^{\rm Re}\wedge \left< \boldsymbol{l} \right>_{\mathbb{R}}\cong V_q^{\rm Re}$. 
\end{proof}

\subsection{The Schouten bracket on $S^4$}\label{s4.5}
The Schouten bracket of $\mathbb{C}^{4}\backslash\{0\}$ satisfies the relation 
\begin{equation*}\label{s4.5e1}
[\Phi(w), \Phi(w^{\prime})]_{Sch}=\Phi[w, w^{\prime}]_{Sch}
\end{equation*}
for any 2-vectors $w, w^{\prime}$ on $\mathbb{C}^{4}\backslash\{0\}$ since $\Phi$ is isometric and preserves the Riemannian connection on $\mathbb{C}^{4} \backslash\{0\}$. 
We recall that the map $\Phi$ defines the real structure on the space $\widetilde{V}_2$ of $\mathbb{C}^{*}$-invariant holomorphic 2-vectors on $\mathbb{C}^{4} \backslash\{0\}$. 
Let $w, w^{\prime}$ be two elements of the real form $\widetilde{V}_2^{\rm Re}$ of $\widetilde{V}_2$. 
Then the Schouten bracket $[w, w^{\prime}]_{Sch}$ is in the real form $\widetilde{V}_3^{\rm Re}$ of $\widetilde{V}_3$ 
since $[w, w^{\prime}]_{Sch}=[\Phi(w), \Phi(w^{\prime})]_{Sch}=\Phi[w, w^{\prime}]_{Sch}$. 
Thus, the Schouten bracket maps $\widetilde{V}_2^{\rm Re}\times \widetilde{V}_2^{\rm Re}$ to $\widetilde{V}_3^{\rm Re}$. 
%\begin{prop}
%The equation $[dp(u),dp(v)]_{Sch}=0$ is equal to $[d\pi(u^{2,0}),d\pi(v^{2,0})]_{Sch}=0$ for $u,v\in \widetilde{W}_2$. 
%\end{prop}
\begin{prop}\label{s4.5p1}
Let $u,v$ be elements of $\widetilde{W}_2$. 
The equation $[dp(u),dp(v)]_{Sch}=0$ is equivalent to $[d\pi(u^{2,0}),d\pi(v^{2,0})]_{Sch}=0$. 
\end{prop}
\begin{proof}
If we take two elements $u,v$ of $\widetilde{W}_2$, 
then the $(3,0)$-part $[u,v]_{Sch}^{3,0}=[u^{2,0},v^{2,0}]_{Sch}$ of the Schouten bracket $[u,v]_{Sch}$ is in $\widetilde{V}_3^{\rm Re}$. 
It follows from Proposition~\ref{s4.4p1} that $[u, v]_{Sch}$ is in $\widetilde{W}_3$. 
The equation $dp([u, v]_{Sch})=0$ is equal to $d\pi([u, v]^{3,0}_{Sch})=0$ 
by the isomorphism $W_3\cong V_3^{\rm Re}$ in Corollary~\ref{s4.4c1}.  
The projection $p$ and $\pi$ preserve the the Schouten bracket $[\cdot,\cdot]_{Sch}$. 
Hence $[dp(u),dp(v)]_{Sch}=dp([u, v]_{Sch})=0$ if and only if $[d\pi(u^{2,0}),d\pi(v^{2,0})]_{Sch}=d\pi([u, v]^{3,0}_{Sch})=0$. 
\end{proof}
%The Schouten bracket induces the map from $\widetilde{W}_2\times \widetilde{W}_2$ to $\widetilde{W}_3$. 
%Moreover, since $dp[u,v]_{Sch}=[dp(u),dp(v)]_{Sch}$ for $u,v\in \widetilde{W}_2$, 
%it induces the map from $W_2\times W_2$ to $W_3$, which is the Schouten bracket of $\mathbb{HP}^1$ restricted to $W_2$. 
%In the previous subsection, we can see that the space $W_2$ is the real form of the space $V_2$ of holomorphic 2-vectors on $\mathbb{CP}^{3}$. 
Proposition~\ref{s4.5p1} implies that any Poisson structure on $\mathbb{HP}^1$ which is in $W_2$ corresponds to 
a holomorphic Poisson structure on $\mathbb{CP}^{3}$ which is in the real form $V_2^{\rm Re}$ of $V_2$. 
Under the identification of $\widetilde{V}_q$ with $S^q\mathbb{C}^{4}\otimes \wedge^q(\mathbb{C}^{4})^{\vee}$, 
the Schouten brackets on $\widetilde{V}_2$ and $V_2$ are given by the maps $\widetilde{F}$ and $F$ as in \S\ref{s3.3}, respectively. 
Then the Schouten bracket on $W_2$ corresponds to the map $F$ restricted to 
$S^2\otimes \wedge^2(gl(2, \mathbb{H})\otimes gl(2, \mathbb{H}))\, / \, S^2\otimes \wedge^2(gl(2, \mathbb{H})\otimes \left< E \right>_{\mathbb{R}})$. 
Hence, we have 
\begin{thm}\label{s4.5t1}
The space of Poisson structures on $\mathbb{HP}^1$ induced by holomorphic Poisson structures on $\mathbb{CP}^{3}$ is identified with 
\[
\{a\in S^2\otimes \wedge^2(gl(2, \mathbb{H})\otimes gl(2, \mathbb{H}))\, /\, S^2\otimes \wedge^2(gl(2, \mathbb{H})\otimes \left< E \right>_{\mathbb{R}}) \mid F(a,a)=0\},
\]
and it is the real form of the space of holomorphic Poisson structures on $\mathbb{CP}^3$. $\hfill\Box$
\end{thm}
The space has the $\mathbb{R}^*$-action induced by the $\mathbb{C}^*$-action on the space of holomorphic Poisson structures on $\mathbb{CP}^{3}$. 
The quadratic polynomial $F(a,a)$ has real coefficients. 
Hence, we obtain 
\begin{cor}\label{s4.5c1}
The quotient of the space of Poisson structures on $\mathbb{HP}^1$ induced by holomorphic Poisson structures on $\mathbb{CP}^{3}$ by the $\mathbb{R}^*$-action 
is a real projective algebraic variety. $\hfill\Box$
\end{cor}

\subsection{Examples of Poisson structures on $S^4$}\label{s4.7}
In order to check that $dp(w)$ is a non-zero 2-vector on $\mathbb{HP}^1$, 
we represent 1-vectors by using the local coordinate of $\mathbb{HP}^1$. 
Let $U'_m$  be the open set $ \{[h_0,h_1] \in  \mathbb{HP}^1 | h_m \ne 0 \}$ of $\mathbb{HP}^1$ for $m=0,1$. 
The functions $s_0 = h_1h_0^{-1}$ and $s_1= h_0h_{1}^{-1}$ provide coordinate on $U'_0$ and $U'_1$, respectively. 
%Then the system $ \{s_0,s_1 \} \backslash  \{s_m \}$  without $s_m$ defines the inhomogeneous coordinate on $U'_m$. 
We can write $h_0=(z_0,z_1)$ and $h_1=(z_2, z_3)$ by using the complex coordinate $(z_0, z_1, z_2, z_3)$ of $\mathbb{C}^4 \backslash\{0\}$. 
%Let $v$ be the real part $a+\bar{a}$ of a holomorphic 2-vector $a=\sum_{i,j,k,l}a_{ijkl}z_iz_j\frac{\partial}{\partial z_k}\wedge \frac{\partial}{\partial z_l}$ such that $(a_{ijkl})$ is in $W$. 
%Then $v$ defines the real 2-vector on $\mathbb{HP}^1$ by the push forward $dp(v)$ by the projection $p :\mathbb{C}^4 \backslash\{0\} \to \mathbb{HP}^1$. 
%In order to compute the Schouten bracket $[dp(v), dp(v)]_{Sch}$ of a 2-vector on $\mathbb{HP}^1$, 
%we represent $dp(v)$ by 
%\[
%dp(v) = dp(a)+dp(\bar{a})
%\]
%on the local coordinate $U'_m$. 
On the local coordinate $(U'_0, s_0)$, if we denote the inhomogeneous coordinate $s_0$ by $t_{0} +jt_{1}$, 
then we have 
\[
dp( \frac{ \partial }{ \partial z_0})=-t_0\frac{\partial}{\partial t_0}-t_1\frac{\partial}{\partial t_1},\ 
dp( \frac{ \partial }{ \partial z_1})=\bar{t}_1\frac{\partial}{\partial t_0}-\bar{t}_0\frac{\partial}{\partial t_1},\ 
dp( \frac{ \partial }{ \partial z_2})=\frac{\partial}{\partial t_0},\ 
dp( \frac{ \partial }{ \partial z_3})=\frac{\partial}{\partial t_1}
\]
%\begin{eqnarray*}
%dp( \frac{ \partial }{ \partial z_0})=-t_0\frac{\partial}{\partial t_0}-t_1\frac{\partial}{\partial t_1},\ 
%dp( \frac{ \partial }{ \partial z_1})=\bar{t}_1\frac{\partial}{\partial t_0}-\bar{t}_0\frac{\partial}{\partial t_1}, \\
%dp( \frac{ \partial }{ \partial z_2})=\frac{\partial}{\partial t_0},\ 
%dp( \frac{ \partial }{ \partial z_3})=\frac{\partial}{\partial t_1}
%\end{eqnarray*}
at $z_0=1, z_1=0, z_2=t_0, z_3=t_1$. 
Similarly, on $U'_1$ with the coordinate $s_1=t_{0} +jt_{1}$, 
we also have 
\[
dp( \frac{ \partial }{ \partial z_0})=\frac{\partial}{\partial t_0},\ 
dp( \frac{ \partial }{ \partial z_1})=\frac{\partial}{\partial t_1}, \ 
dp( \frac{ \partial }{ \partial z_2})=-t_0\frac{\partial}{\partial t_0}-t_1\frac{\partial}{\partial t_1}, \ 
dp( \frac{ \partial }{ \partial z_3})=\bar{t}_1\frac{\partial}{\partial t_0}-\bar{t}_0\frac{\partial}{\partial t_1}
\]
%\begin{eqnarray*}
%dp( \frac{ \partial }{ \partial z_0})&=&\frac{\partial}{\partial t_0}, \\
%dp( \frac{ \partial }{ \partial z_1})&=&\frac{\partial}{\partial t_1}, \\
%dp( \frac{ \partial }{ \partial z_2})&=&-t_0\frac{\partial}{\partial t_0}-t_1\frac{\partial}{\partial t_1}, \\
%dp( \frac{ \partial }{ \partial z_3})&=&\bar{t}_1\frac{\partial}{\partial t_0}-\bar{t}_0\frac{\partial}{\partial t_1}
%\end{eqnarray*}
at $z_0=t_0, z_1=t_1, z_2=1, z_3=0$. 
%Thus we can only take $\frac{ \partial}{\partial t_0} \wedge \frac{\partial}{ \partial t_1} \wedge  \frac{\partial}{ \partial \bar{t_l}}$ for $l=0,1$ and the complex conjugate of it on each $W_m$.  
%Hence $[dp(v),dp(v)]_{Sch}=[dp(a),dp(\bar{a})]_{Sch}+[dp(\bar{a}),dp(a)]_{Sch}.$  It vanishes if and only if $[dp(a),dp(\bar{a})]_{Sch}$ is equal to zero for each $m$.

Let $u$ and $v$ be two complex 1-vectors on $\mathbb{C}^4\backslash \{0\}$. 
If $u+\bar{u}$ and $v+\bar{v}$ are $\mathbb{H}^*$-invariant, 
then the Schouten bracket of the real 2-vector $dp(u+\bar{u})\wedge dp(v+\bar{v})$ on $\mathbb{HP}^1$ is given by 
$2dp([u, v]+\overline{[u, v]})\wedge (u+\bar{u})\wedge (v+\bar{v}))$. 
%\[
%2dp([u+\bar{u}, v+\bar{v}]\wedge (u+\bar{u})\wedge (v+\bar{v}))=2dp([u, v]+\bar{[u, v]})\wedge (u+\bar{u})\wedge (v+\bar{v})).
%\]
%\begin{eqnarray*}
%&&[dp(u+\bar{u})\wedge dp(v+\bar{v}),dp(u+\bar{u})\wedge dp(v+\bar{v})]_{Sch}\\
%&=&2dp([u+\bar{u}, v+\bar{v}]\wedge (u+\bar{u})\wedge (v+\bar{v}))\\
%&=&2dp([u, v]+\bar{[u, v]})\wedge (u+\bar{u})\wedge (v+\bar{v})).
%\end{eqnarray*}
It implies that $dp(u+\bar{u})\wedge dp(v+\bar{v})$ defines a Poisson structure in the case of $[u,v]=0$. 
In general, if we take $n$ 1-vectors $v_1, \dots , v_n$ on $\mathbb{C}^4\backslash \{0\}$ such that $[v_i, v_j]=0$ for $i,j=1, \dots,n$, 
then a real 2-vector $\sum_{i<j} a_{ij}dp(v_i+\bar{v}_i)\wedge dp(v_j+\bar{v}_j)$ is a Poisson structure 
for a constant $a_{ij}\in \mathbb{R}$ by the similar argument of $\mathbb{CP}^3$. 

\begin{ex}\label{s4.7ex1}
{\rm 
We define 1-vectors $v_1, v_2, v_3$ by 
\[
v_1=z_0 \frac{\partial}{\partial z_0}+z_1 \frac{\partial}{\partial z_1},\ 
v_2=i(z_0 \frac{\partial}{\partial z_0}-z_1 \frac{\partial}{\partial z_1}),\ 
v_3=i(z_2 \frac{\partial}{\partial z_2}-z_3 \frac{\partial}{\partial z_3})
\]
on $\mathbb{C}^4\backslash \{0\}$. Then $[v_1, v_2]=[v_1, v_3]=[v_2, v_3]=0$. 
We can also check that each push forward $dp((v_1+\bar{v}_1)\wedge (v_2+\bar{v}_2)), dp((v_1+\bar{v}_1)\wedge (v_3+\bar{v}_3)), dp((v_2+\bar{v}_2)\wedge (v_3+\bar{v}_3))$ 
is a non-zero 2-vector on $\mathbb{HP}^1$ by using the local coordinate $U'_0$. 
Hence we obtain a family of Poisson structures 
\[
a_1dp((v_1+\bar{v}_1)\wedge (v_2+\bar{v}_2))+a_2 dp((v_1+\bar{v}_1)\wedge (v_3+\bar{v}_3))+a_3 dp((v_2+\bar{v}_2)\wedge (v_3+\bar{v}_3))
\]
for $a_1, a_2, a_3 \in \mathbb{R}$. 
}
\end{ex}

\section{Poisson structures on $\mathbb{CP}^m$ and $\mathbb{HP}^n$}
In this section, we generalize the above argument to the higher dimensional case $\mathbb{CP}^m$ and $\mathbb{HP}^n$. 
\subsection{Holomorphic Poisson structures on $\mathbb{CP}^m$}
We denote by $(z_0, z_1,\dots, z_m)$ the standard coordinate of $\mathbb{C}^{m+1}$. 
Let $\pi$ be the projection $\mathbb{C}^{m+1}\backslash\{0\}\to \mathbb{CP}^m$. 
The action of $\mathbb{C}^*$ on $\mathbb{C}^{m+1}\backslash\{0\}$ induces the vector field 
\[
\boldsymbol{l}=z_0\frac{\partial}{\partial z_0}+z_1\frac{\partial}{\partial z_1}+\dots+z_m\frac{\partial}{\partial z_m} 
\]
on $\mathbb{C}^{m+1}\backslash\{0\}$. 
We can regard the space of holomorphic 2-vectors on $ \mathbb{CP}^m$ as the quotient space 
\[
\Bigl\{\sum_{i,j,k,l}a_{ijkl}z_iz_j\frac{\partial}{\partial z_k}\wedge \frac{\partial}{\partial z_l} \Bigm| a_{ijkl}\in \mathbb{C}\Bigr\} \Bigm/ 
\Bigl\{(\sum_{i,j}a_{ij}z_i\frac{\partial}{\partial z_j})\wedge \boldsymbol{l} \Bigm| a_{ij}\in \mathbb{C}\Bigr\}
\]
%of the space $\{\sum_{i,j,k,l}a_{ijkl}z_iz_j\frac{\partial}{\partial z_k}\wedge \frac{\partial}{\partial z_l} \mid a_{ijkl}\in \mathbb{C}\}$ divided by 
%the subspace $\{(\sum_{i,j}a_{ij}z_i\frac{\partial}{\partial z_j})\wedge(\sum_{k}z_k\frac{\partial}{\partial z_k})\mid a_{ij}\in \mathbb{C}\}$. 
which is identified with 
\[
S^2\mathbb{C}^{m+1}\otimes \wedge^2(\mathbb{C}^{m+1})^{\vee}\, / \, S^2\otimes \wedge^2(sl(m+1,\mathbb{C})\otimes \left< E \right>_{\mathbb{C}})
\]
where $S^2\otimes \wedge^2$ is the projection 
from $\displaystyle gl(\mathbb{C}^{m+1})\otimes gl(\mathbb{C}^{m+1})\cong (\otimes^2\mathbb{C}^{m+1})\otimes (\otimes^2 (\mathbb{C}^{m+1})^{\vee})$ 
to $S^2\mathbb{C}^{m+1}\otimes \wedge^2(\mathbb{C}^{m+1})^{\vee}$. 
The Schouten bracket for $\mathbb{C}^*$-invariant 2-vectors on $\mathbb{C}^{m+1}\backslash\{0\}$ can be considered as the map 
\[
\widetilde{F} : S^2\mathbb{C}^{m+1}\otimes \wedge^2(\mathbb{C}^{m+1})^{\vee}\times S^2\mathbb{C}^{m+1}\otimes \wedge^2(\mathbb{C}^{m+1})^{\vee}
\to S^3\mathbb{C}^{m+1}\otimes \wedge^3(\mathbb{C}^{m+1})^{\vee}. 
\]
The Schouten bracket on $\mathbb{CP}^m$ is given by 
the map from 
$S^2\mathbb{C}^{m+1}\otimes \wedge^2(\mathbb{C}^{m+1})^{\vee} / \{S^2\otimes \wedge^2(sl(m+1,\mathbb{C})\otimes \left< E \right>_{\mathbb{C}})\} 
\times S^2\mathbb{C}^{m+1}\otimes \wedge^2(\mathbb{C}^{m+1})^{\vee} / \{S^2\otimes \wedge^2(sl(m+1,\mathbb{C})\otimes \left< E \right>_{\mathbb{C}})\}$ 
to $S^3\mathbb{C}^{m+1}\otimes \wedge^3(\mathbb{C}^{m+1})^{\vee} / \{S^3\otimes \wedge^3(\otimes^{2}sl(m+1,\mathbb{C})\otimes \left< E \right>_{\mathbb{C}})\}$ induced by $\widetilde{F}$. 
We denote by $F$ the induced map. 
Then, we have 
\begin{prop}\label{s5.1p1}
The space of holomorphic Poisson structures on $\mathbb{CP}^m$ is identified with 
\[
\{a\in S^2\mathbb{C}^{m+1}\otimes \wedge^2(\mathbb{C}^{m+1})^{\vee}\, / \, S^2\otimes \wedge^2(sl(m+1,\mathbb{C})\otimes \left< E \right>_{\mathbb{C}}) \mid F(a,a)=0\}.
\]
Moreover, the quotient of the space by the $\mathbb{C}^*$-action is a projective algebraic variety.
$\hfill\Box$
\end{prop}

%\subsubsection{Examples of Poisson structures on $\mathbb{CP}^m$}
We provide some examples of Poisson structures on $\mathbb{CP}^m$. 
\begin{ex}\label{s5.1ex1}
{\rm 
 We define $v_1$ by the 2-vector 
\[
v_1=z_0z_i \frac{\partial }{\partial z_j}\wedge \frac{\partial }{\partial z_k}
\]
on $\mathbb{C}^{m+1}\backslash \{0\}$ for $1\le i,j,k\le m$. 
Then the Schouten bracket of $d\pi(v_1)$ vanishes on $\mathbb{CP}^m$. 
Hence the push forward $d\pi(v_1)$ is a Poisson structure on $\mathbb{CP}^m$. 
}
\end{ex}

\begin{ex}\label{s5.1ex2}
{\rm 
We take the 2-vector 
\[
v_2=(z_0 \frac{\partial }{\partial z_1}- z_1\frac{\partial }{\partial z_0}) \wedge 
(z_i \frac{\partial }{\partial z_j}- z_j\frac{ \partial }{ \partial z_i})
\]
on $\mathbb{C}^{m+1}\backslash \{0\}$ for $1\le i<j\le m$. 
Then $d\pi(v_2)$ is a Poisson structure on $\mathbb{CP}^m$. 
}
\end{ex}

\begin{ex}\label{s5.1ex3}
{\rm 
We define $v_3$ by the 2-vector 
\[
v_3=(z_0\frac{\partial }{\partial z_1}+z_1\frac{\partial }{\partial z_0}+z_i\frac{\partial }{\partial z_{i+1}}+z_{i+1}\frac{\partial }{\partial z_i})\wedge 
(z_0\frac{\partial }{\partial z_{i+1}}+z_1\frac{\partial }{\partial z_i}+z_i\frac{\partial }{\partial z_1}+z_{i+1}\frac{\partial }{\partial z_0})
\]
on $\mathbb{C}^{m+1}\backslash \{0\}$ for $2\le i\le m$. 
Then $d\pi(v_3)$ is a Poisson structure on $\mathbb{CP}^m$. 
}
\end{ex}

\begin{ex}\label{s5.1ex4}
{\rm 
We define $u_1, u_2,\dots ,u_m$ by the 1-vectors 
\[
u_1=z_1\frac{\partial }{\partial z_1}\ ,\quad u_2=z_2\frac{\partial }{\partial z_2}\ ,\dots\dots,\ u_m=z_m\frac{\partial }{\partial z_m}
\]
on $\mathbb{C}^{m+1}\backslash \{0\}$. 
Then, we obtain a family of the Poisson structures 
\[
a_{12}d\pi(u_1\wedge u_2)+a_{13} d\pi(u_1\wedge u_3)+\dots+ a_{m-1\, m} d\pi(u_{m-1}\wedge u_m)
\] 
for $a_{12}, a_{13},\dots, a_{m-1\, m} \in \mathbb{C}$. 
}
\end{ex}

\subsection{Poisson structures on $\mathbb{HP}^n$}
The space $\mathbb{H}^{n+1}$ is isomorphic to the right $\mathbb{C}$-module $\mathbb{C}^{2n+2}$ 
under the identification of $\mathbb{H}$ with $\mathbb{C}^2 (\cong \mathbb{C} + j \mathbb{C})$. 
The right action of $j$ on $\mathbb{C}^{2n+2} \backslash\{0\}$ induces the map $\boldsymbol{j} : \mathbb{C}^{2n+2}\backslash\{0\} \to \mathbb{C}^{2n+2}\backslash \{0\}$ given by 
\[
(z_0,z_1,\dots,z_{2n},z_{2n+1}) \mapsto (- \bar{z}_1,\bar{z}_0,\dots,-\bar{z}_{2n+1},\bar{z}_{2n})
\]
for $(z_0,z_1,\dots,z_{2n},z_{2n+1})\in \mathbb{C}^{2n+2}\backslash\{0\}$. 
Let $p$ be the projection from $\mathbb{C}^{2n+2} \backslash\{0\}\cong \mathbb{H}^{n+1}\backslash\{0\}$ to $\mathbb{HP}^n$. 
We consider real $q$-vectors on $\mathbb{HP}^n$ given by 
the push forward of $\mathbb{H}^*$-invariant real $q$-vectors with holomorphic $(q,0)$-parts on $\mathbb{C}^{2n+2}\backslash\{0\}$. 

We define the map $\Phi: \Gamma(T^{\mathbb{C}}\mathbb{C}^{2n+2}\backslash\{0\}) \to \Gamma(T^{\mathbb{C}}\mathbb{C}^{2n+2}\backslash\{0\})$ by 
\[
\Phi(u)=\overline{d\boldsymbol{j}(u)}
\]
for any $u\in \Gamma(T^{\mathbb{C}}\mathbb{C}^{2n+2}\backslash\{0\})$. 
Then the map $\Phi$ induces a real structure on the space of the holomorphic $2$-vectors on $\mathbb{CP}^{2n+1}$. 
We obtain 
\begin{thm}\label{s5.2t1}
The space of Poisson structures on $\mathbb{HP}^n$ induced by holomorphic Poisson structures on $\mathbb{CP}^{2n+1}$ is identified with 
\[
\{a\in S^2\otimes \wedge^2(gl(n+1, \mathbb{H})\otimes gl(n+1, \mathbb{H}))\, /\, S^2\otimes \wedge^2(gl(n+1, \mathbb{H})\otimes \left< E \right>_{\mathbb{R}}) \mid F(a,a)=0\},
\]
and it is the real form of the space of holomorphic Poisson structures on $\mathbb{CP}^{2n+1}$. 
Moreover, the quotient of the space by the $\mathbb{R}^*$-action is a real projective algebraic variety. $\hfill\Box$
\end{thm}

%\subsubsection{Examples of Poisson structures on $\mathbb{HP}^n$}
We provide some examples of Poisson structures on $\mathbb{HP}^n$. 
\begin{ex}\label{s5.2ex1}
{\rm 
We take 1-vectors 
\[
u=i(z_0 \frac{\partial}{\partial z_0}-z_1 \frac{\partial}{\partial z_1}),\ v=i(z_{2k} \frac{\partial}{\partial z_{2k}} -z_{2k+1}\frac{\partial}{\partial z_{2k+1}})
\]
on $\mathbb{C}^{2n+2}\backslash \{0\}$ for $k=1,\dots,n$. 
Then $dp(u+\bar{u})\wedge dp(v+\bar{v})$ is a Poisson structure on $\mathbb{HP}^n$. 
}
\end{ex}

\begin{ex}\label{s5.2ex2}
{\rm 
We define 1-vectors $v_1,\dots, v_{n+1}$ by 
\[
v_1=z_0 \frac{\partial}{\partial z_0}+z_1 \frac{\partial}{\partial z_1},\ 
v_2=i(z_0 \frac{\partial}{\partial z_0}-z_1 \frac{\partial}{\partial z_1}),\ \dots, \ 
v_{n+1}=i(z_{2n} \frac{\partial}{\partial z_{2n}}-z_{2n+1} \frac{\partial}{\partial z_{2n+1}})
\]
on $\mathbb{C}^{2n+2}\backslash \{0\}$. 
We obtain a family of Poisson structures 
\begin{multline*}
a_{12}dp((v_1+\bar{v}_1)\wedge (v_2+\bar{v}_2))+a_{13}dp((v_1+\bar{v}_1)\wedge (v_3+\bar{v}_3))+\cdots \\
\cdots +a_{n-1\, n+1}dp((v_{n-1}+\bar{v}_{n-1})\wedge (v_{n+1}+\bar{v}_{n+1}))+a_{n\, n+1}dp((v_{n}+\bar{v}_{n})\wedge (v_{n+1}+\bar{v}_{n+1}))
\end{multline*}
for $a_{12},a_{13}, \dots, a_{n-1\, n+1}, a_{n\, n+1} \in \mathbb{R}$. 
}
\end{ex}

\section{Poisson structures on $S^4$ induced by holomorphic foliations on $\mathbb{CP}^3$}
In this section, we provide examples of Poisson structures on $S^4$ 
associated with Cerveau and Lins Neto's classification of holomorphic foliations. 

\subsection{Holomorphic Poisson structures on $\mathbb{CP}^3$ induced by holomorphic foliations}
There exist rich holomorphic foliations of codimension $1$ on $\mathbb{CP}^3$. 
Especially, holomorphic foliations of degree $2$ are completely classified and the space of such foliations is decomposed into the six components~\cite{CLN}. 
Such a foliation is given by a holomorphic 1-form $\alpha$ on $\mathbb{C}^4\backslash \{0\}$. 
The $1$-form $\alpha$ induces a 2-vector $w$ by the correspondence 
\begin{equation}
\alpha=i_{\boldsymbol{l}\wedge w}(dz_0\wedge dz_1\wedge dz_2\wedge dz_3)\label{s6e0}
\end{equation}
on $\mathbb{C}^4\backslash \{0\}$. 
Such a 2-vector $w$ is not uniquely determined by the $1$-form $\alpha$ 
since there exists an ambiguity for 2-vectors factorized by $\boldsymbol{l}$. 
The correspondence $w \mapsto \alpha$ induces the map $S^2\mathbb{C}^{4}\otimes \wedge^2(\mathbb{C}^{4})^{\vee} \to S^3\mathbb{C}^{4}\otimes (\mathbb{C}^{4})^{\vee}$. 
It is well-known that the push forward $d\pi(w)$ is the Poisson structure on $\mathbb{CP}^3$. 
Actually, the integrability condition $\alpha\wedge d\alpha =0$ is equivalent to $[w,w]_{Sch}\wedge \boldsymbol{l}=0$. 
%Hence, we can construct many Poisson structures corresponding to holomorphic foliations on $\mathbb{CP}^3$. 
Pym provided a correspondence of the $1$-form $\alpha$ to a {\it unimoduler} Poisson bracket $\{,\}$ on $\mathbb{C}^{4}$ by 
\[
\{f,g\}=\frac{df\wedge dg \wedge d\alpha}{dz_0\wedge dz_1\wedge dz_2\wedge dz_3}
\]
for functions $f, g$~\cite{P}. The Poisson bracket $\{,\}$ induces a Poisson structure $d\pi(w)$ on $\mathbb{CP}^3$ given by 
\begin{equation}
w=\sum_{i,j=0}^{3}\frac{dz_i \wedge dz_j \wedge d\alpha}{dz_0\wedge dz_1\wedge dz_2\wedge dz_3}\frac{\partial }{\partial z_i}\wedge \frac{\partial }{\partial z_j}\label{s6e1}
\end{equation}
on $\mathbb{C}^{4}\backslash \{0\}$. 

\subsection{Poisson structures on $S^4$ induced by holomorphic foliations on $\mathbb{CP}^3$}
We define a transformation $\Phi$ on the space of $p$-forms on $\mathbb{C}^4\backslash \{0\}$ by 
\[
\Phi(\theta)=\overline{\boldsymbol{j}^*\theta}
\]
for a $p$-form $\theta$. 
We consider a holomorphic 1-form on $\mathbb{C}^4\backslash \{0\}$ 
associated with a holomorphic foliation of degree $2$ on $\mathbb{CP}^3$. 
The coefficients of such a 1-form are polynomials of degree $3$. 
The transformation $\Phi$ is a real structure on the space of 1-forms whose coefficients are polynomials of degree $3$ on $\mathbb{C}^4\backslash \{0\}$. 
Let $\alpha$ be a 1-form on $\mathbb{C}^4\backslash \{0\}$ 
associated with a holomorphic foliation of degree $2$ on $\mathbb{CP}^3$. 
The 1-form $\alpha$ induces a 2-vector $w$ such that $d\pi(w)$ is the Poisson structure on $\mathbb{CP}^3$ by the correspondence (\ref{s6e0}). 
If the 2-vector $w$ is real with respect to the real structure defined in \S\ref{s4.4}, then the 1-form $\alpha$ is real with respect $\Phi$. 
In general, the converse is not true by the choice of the 2-vector $w$ for the 1-form $\alpha$. 
However, under the correspondence (\ref{s6e1}) it is true : 
\begin{thm}\label{s6t1}
The $2$-vector $w$ given by (\ref{s6e1}) is real if the $1$-from $\alpha$ is real. 
\end{thm}
\begin{proof} 
It follows from $\Phi(dz_i)=(-1)^{i+1}dz_{\varphi(i)}$ and $\Phi(\frac{\partial }{\partial z_i})=(-1)^{i}\frac{\partial }{\partial z_{\varphi(i)}}$ 
that $\Phi(w)=w$ if $\Phi(\alpha)=\alpha$. 
\end{proof}
Therefore, any 1-form $\alpha$ which is real with respect to $\Phi$ associated with a foliation of degree $2$ on $\mathbb{CP}^3$ induces a Poisson structure on $S^4$.

\subsection{Examples of Poisson structures on $S^4$ associated with holomorphic foliations on $\mathbb{CP}^3$}\label{s6.3}
The space of holomorphic foliations of degree $2$ on $\mathbb{CP}^3$ is decomposed into 
the six components $L(1,1,1,1), L(1,1,2), R(2,2), R(1,3),S(2,3), E(3)$ 
%$\overline{L(1,1,1,1)}, \overline{L(1,1,2)}, \overline{R(2,2)}, \overline{R(1,3)},S(2,3), \overline{E(3)}$ 
(we refer to \cite{CLN} for the definitions of these components). 
We provide a 1-form which is real with respect to $\Phi$ associated with such components. 

\vspace{\baselineskip}
\noindent
\textbf{The $L(1,1,1,1)$ component} 

\vspace{0.2\baselineskip}
\noindent
The component $L(1,1,1,1)$ is given by the closure of the ${\rm GL}(4, \mathbb{C})$-orbit of the following one form 
\[
\alpha=z_0 z_1 z_2 z_3\sum_{i=0}^3 a_i \frac{dz_i}{z_i}
\]
for $a_0,\dots, a_3 \in \mathbb{C}$ such that $\sum_{i=0}^3 a_i=0$. 
The form $\alpha$ is real with respect to $\Phi$ if and only if $a_1=\overline{a}_0, a_3=\overline{a}_2$ and ${\rm Re}\, a_2=-{\rm Re}\, a_0$, that is, 
\[
\alpha=z_0 z_1 z_2 z_3\left( a_0 \frac{dz_0}{z_0}+\overline{a}_0 \frac{dz_1}{z_1}+a_2 \frac{dz_2}{z_2}+\overline{a}_2 \frac{dz_3}{z_3}\right)
\]
for $a_0, a_2\in \mathbb{C}$ with ${\rm Re}\, a_2=-{\rm Re}\, a_0$. 
Then the 2-vector $w$ given by (\ref{s6e1}) induces the Poisson structure on $S^4$. 

\vspace{\baselineskip}
\noindent
\textbf{The $L(1,1,2)$ component} 

\vspace{0.2\baselineskip}
\noindent 
We identify $S^i\mathbb{C}^4$ with the space of homogeneous polynomials of degree $i$ on $\mathbb{C}^4\backslash \{0\}$. 
We define an one form $\alpha$ by 
\[
\alpha=z_0 z_1 f\left( a_0 \frac{dz_0}{z_0}+a_1 \frac{dz_1}{z_1}+a_2 \frac{df}{f}\right)
\]
for $f\in S^2\mathbb{C}^4$ and $a_0,a_1, a_2 \in \mathbb{C}$ such that $a_0+a_1+2 a_2=0$. 
The one form $\alpha$ induces the foliation of the $L(1,1,2)$ component. 
The form $\alpha$ is real if and only if 
\begin{equation}
\left\{
\begin{array}{ll}
f=cz_0z_1,& |a_0|\neq |a_1|, \\
f=f_2+f_1z_0+\frac{\overline{a}_1}{a_0}\Phi(f_1)z_1+\frac{\sqrt{-1}\,r}{a_0-a_1},& |a_0|=|a_1|,\ a_0\neq a_1,  \\ 
f=f_2+f_1z_0+\frac{\overline{a}_1}{a_0}\Phi(f_1)z_1+c',& a_0= a_1  
\end{array}
\right.
\end{equation}
for $r\in\mathbb{R}$, $c,c'\in \mathbb{C}$, $f_1\in S^1 \left<z_2,z_3\right>_{\mathbb{C}}$ and $f_2\in S^2\left<z_2,z_3\right>_{\mathbb{C}}$ 
with $\Phi(f_2)=-\frac{\overline{a}_1}{a_0}f_2$, where $S^i\left<z_2,z_3\right>_{\mathbb{C}}$ is the space of homogeneous polynomials of degree $i$ of $z_2$ and $z_3$.

\vspace{\baselineskip}
\noindent
\textbf{The $R(2,2)$ component} 

\vspace{0.2\baselineskip}
\noindent
The component $R(2,2)$ is consisted by 1-forms 
\[
\alpha=gdf-fdg
\]
for $f,g\in S^2\mathbb{C}^4$. 
If we take $f,g$ in $gl(2, \mathbb{H})\cap S^2\mathbb{C}^4$, 
then $\alpha=gdf-fdg$ is real with respect to $\Phi$. 
%and the corresponding 2-vector $w$ defines the Poisson structure on $S^4$. 
For example, $(f,g)=(z_0^2+z_1^2, \sqrt{-1}\,z_0z_1), (z_0^2+z_1^2+z_2^2+z_3^2, \sqrt{-1}\,z_0z_1), (z_0^2+z_1^2+z_2^2+z_3^2, \sqrt{-1}\,(z_0z_1+z_2z_3))$. \\
%induce Poisson structures on $S^4$. 

\vspace{\baselineskip}
\noindent
\textbf{The $R(1,3)$ component} 

\vspace{0.2\baselineskip}
\noindent
Let $\alpha$ be a 1-form 
\[
\alpha=3fdz_0-z_0df
\]
for $f\in S^3\mathbb{C}^4$. 
Then the form $\alpha$ is in the component $R(1,3)$. 
The form $\alpha$ is real if and only if $f=c'z_0^3+cz_0^2z_1+\sqrt{-1}\,rz_0z_1^2+\overline{c}z_1^3$ for $r\in\mathbb{R}$ and $c,c'\in \mathbb{C}$, i.e., 
\[
\alpha=(cz_0^2z_1+irz_0z_1^2+\overline{c}z_1^3)dz_0-(cz_0^3+\sqrt{-1}\,rz_0^2z_1+\overline{c}z_0z_1^2)dz_1
\]
for $r\in\mathbb{R}$ and $c\in \mathbb{C}$. The example is also contained in the $S(2,3)$ component.

\vspace{\baselineskip}
\noindent
\textbf{The $S(2,3)$ component}

\vspace{0.2\baselineskip}
\noindent 
We define a 1-form $\alpha$ by 
\[
\alpha=f_0dz_0+f_1dz_1+f_2dz_2
\]
for $f_0,f_1,f_2\in S^3 \left<z_0,z_1,z_2\right>_{\mathbb{C}}$. 
%for homogeneous polynomials $f_0,f_1,f_2$ of degree $3$ of $z_0, z_1$ and $z_2$. 
%The form $\alpha$ is in the component $S(2,3)$. 
The form $\alpha$ is real if and only if $f_2=0$ and $f_1=-\Phi(f_0)$, i.e., 
\[
\alpha=f_0dz_0-\Phi(f_0)dz_1
\]
for $f_0\in S^3 \left<z_0,z_1\right>_{\mathbb{C}}$. 

\begin{rem}
{\rm 
The set of Poisson structures corresponding to the exceptional component $E(3)$ is given by 
the closure of the ${\rm GL}(4, \mathbb{C})$-orbit of the Poisson structure 
$(z_0\frac{\partial }{\partial z_1}+z_1\frac{\partial }{\partial z_2}+z_2\frac{\partial }{\partial z_3})\wedge 
(5z_0\frac{\partial }{\partial z_0}+z_1\frac{\partial }{\partial z_1}-3z_2\frac{\partial }{\partial z_2}-7z_3\frac{\partial }{\partial z_3})$~\cite{P}. 
%$X\wedge Y$ for 
%\[
%X=-\frac{5}{4}z_0\frac{\partial }{\partial z_0}-\frac{1}{4}z_1\frac{\partial }{\partial z_1}+\frac{3}{4}z_2\frac{\partial }{\partial z_2}+\frac{7}{4}z_3\frac{\partial }{\partial z_3},\ 
%Y=4z_0\frac{\partial }{\partial z_1}+4z_1\frac{\partial }{\partial z_2}+4z_2\frac{\partial }{\partial z_3}.
%\]
The Poisson structure is not real. 
}
\end{rem}

\noindent
\textbf{Acknowledgements}. 
The authors would like to thank Professor T. Taniguchi for his useful comments. 
They also would like to thank Professor R. Goto for his comment about the relationship 
between holomorphic foliations and holomorphic Poisson structures on $\mathbb{CP}^3$. 
The first named author is supported by Grant-in-Aid for Young Scientists (B) $\sharp$17K14187 from JSPS.

\begin{center}

\end{center}

\vspace{\baselineskip}
\begin{flushright}
\begin{tabular}{l}
\textsc{Takayuki Moriyama}\\
Department of Mathematics\\
Mie University\\
Mie 514-8507, Japan\\
E-mail: takayuki@edu.mie-u.ac.jp\\
\\

\textsc{Takashi Nitta}\\
Department of Mathematics\\
Mie University\\
Mie 514-8507, Japan\\
E-mail: nitta@edu.mie-u.ac.jp
\end{tabular}
\end{flushright}

\end{document}